\newif\ifjsc
\jscfalse

\ifjsc
  \RequirePackage{fix-cm}
  \documentclass[smallextended,envcountsame,envcountsect]{svjour3}
  \smartqed
\else
  \documentclass[11pt]{article}
  \usepackage[margin=1in]{geometry}
\fi

\usepackage{amsmath}
\usepackage{amssymb}
\ifjsc\else\usepackage{amsthm}\fi
\usepackage{mathtools}
\usepackage{array}
\usepackage{booktabs}
\usepackage{graphicx}
\usepackage{placeins}
\usepackage{etoolbox}
\usepackage[numbers,sort&compress]{natbib}
\usepackage{hyperref}
\usepackage[capitalise,noabbrev]{cleveref}
\allowdisplaybreaks
\numberwithin{equation}{section}
\newcolumntype{L}[1]{>{\raggedright\arraybackslash}p{#1}}

\ifjsc
  \spnewtheorem{assumption}[theorem]{Assumption}{\bfseries}{\itshape}
  \AtEndEnvironment{proof}{\qed}
  \renewcommand{\subclassname}{{\bfseries Mathematics Subject Classification (2020)}\enspace}
  \journalname{Journal of Scientific Computing}
  \newcommand{\figwidth}{\textwidth}
  \newcommand{\captionbelow}[1]{\par\vskip\tabcapgap\tabcapgap=0pt\caption{#1}}
\else
  \newtheorem{theorem}{Theorem}[section]
  \newtheorem{lemma}[theorem]{Lemma}
  \newtheorem{proposition}[theorem]{Proposition}
  \newtheorem{corollary}[theorem]{Corollary}
  \newtheorem{assumption}[theorem]{Assumption}
  \theoremstyle{definition}
  \newtheorem{remark}[theorem]{Remark}
  \newtheorem{example}[theorem]{Example}
  \newcommand{\figwidth}{0.8\textwidth}
  \newcommand{\captionbelow}[1]{\caption{#1}}
  \newcommand{\keywords}[1]{\par\medskip\noindent\textbf{Keywords.} {\def\and{\unskip, }#1}}
  \newcommand{\subclass}[1]{\par\smallskip\noindent\textbf{MSC 2020.} {\def\and{\unskip, }#1}}
  \newenvironment{acknowledgements}{\section*{Acknowledgments}}{}
\fi

\newcommand{\R}{\mathbb{R}}
\newcommand{\C}{\mathbb{C}}
\newcommand{\E}{\mathbb{E}}
\newcommand{\cK}{\mathcal{K}}
\newcommand{\cA}{\mathcal{A}}
\newcommand{\cH}{\mathcal{H}}
\newcommand{\cD}{\mathcal{D}}
\newcommand{\cV}{\mathcal{V}}
\newcommand{\cL}{\mathcal{L}}
\newcommand{\cB}{\mathcal{B}}
\newcommand{\cQ}{\mathcal{Q}}
\newcommand{\cS}{\mathcal{S}}
\newcommand{\cZ}{\mathcal{Z}}
\newcommand{\cE}{\mathcal{E}}
\newcommand{\cP}{\mathcal{P}}
\newcommand{\cM}{\mathcal{M}}
\newcommand{\cN}{\mathcal{N}}
\newcommand{\dom}{\operatorname{dom}}
\newcommand{\Span}{\operatorname{span}}
\newcommand{\dt}{\Delta t}

\newcommand{\ip}[2]{\langle #1,#2\rangle}
\newcommand{\Proj}{\Pi}
\newcommand{\Tor}{\mathbb{T}}
\newcommand{\QR}{\mathcal{Q}^{\mathrm{R}}}

\hypersetup{
  pdftitle={A Second-Order Stochastic Dynamic Mode Decomposition by Richardson Extrapolation},
  pdfauthor={Y. Xu}
}

\crefname{theorem}{Theorem}{Theorems}
\crefname{lemma}{Lemma}{Lemmas}
\crefname{proposition}{Proposition}{Propositions}
\crefname{corollary}{Corollary}{Corollaries}
\crefname{remark}{Remark}{Remarks}
\crefname{example}{Example}{Examples}
\crefname{assumption}{Assumption}{Assumptions}

\begin{document}

\title{A Second-Order Stochastic Dynamic Mode Decomposition by Richardson Extrapolation}

\ifjsc
\titlerunning{Second-order SDMD by Richardson extrapolation}
\author{Yuanchao Xu}
\institute{Y. Xu \at Department of Mathematics, Graduate School of Science, Kyoto University, Kitashirakawa Oiwake-cho, Sakyo-ku, Kyoto 606-8502, Japan\\
\email{xu.yuanchao.3a@kyoto-u.ac.jp}}
\date{Received: date / Accepted: date}
\else
\author{Yuanchao Xu\thanks{Department of Mathematics, Graduate School of Science, Kyoto University, Kitashirakawa Oiwake-cho, Sakyo-ku, Kyoto 606-8502, Japan (\texttt{xu.yuanchao.3a@kyoto-u.ac.jp}).}}
\date{}
\fi

\maketitle

\begin{abstract}
Stochastic dynamic mode decomposition (SDMD) approximates the Koopman semigroup of an It\^o diffusion by a Galerkin generator matrix on a dictionary of observables, computed from samples of the invariant law and generator evaluations, and advances observables with a first-order stochastic Taylor step (SDMD-1). The second-order Taylor step (SDMD-2) also needs a third matrix, a dictionary with four derivatives, and a fourth-order differential operator evaluated at every sample. Richardson extrapolation of the SDMD-1 step (SDMD-R) is second-order accurate with the data, the matrices, and the regularity of SDMD-1. The error bound, for exact and empirical matrices, consists of the Galerkin error of the dictionary, a time-stepping term, and a data term of order $m^{-1/2}$ in the sample size $m$. Relative to the Galerkin semigroup, the limit of the schemes as the time step $\dt\to0$, the time-stepping error is $O(\dt^2)$ for SDMD-R and $O(\Theta_N\dt+\dt^2)$ for SDMD-2, with $\Theta_N$ the norm of the second-order Galerkin residual $\Proj_N\cA(I-\Proj_N)\cA|_{\cV_N}$, and the first-order term is attained. For reversible periodic diffusions and trigonometric dictionaries with cutoff $K$, the step restriction $\dt\lesssim K^{-2}$ makes all three schemes contractions, and the error against the Koopman semigroup on $H^\ell$ observables is $O(K^{-(\ell-2)}+\dt^p)$ plus a data term of order $K^{d+2}(m^{-1}\log K)^{1/2}$ in dimension $d$, with $p=1$ for SDMD-1 and $p=2$ for SDMD-R and SDMD-2. In experiments with trigonometric dictionaries the errors of SDMD-R and SDMD-2 differ by at most ten percent. For a spline dictionary with two derivatives, SDMD-R remains second-order accurate, whereas SDMD-2 with pointwise evaluations is first-order accurate.
\subclass{37M10 \and 47D06 \and 60H10 \and 65C30 \and 65P99}
\end{abstract}

\section{Introduction}\label{sec:intro}
Many systems in physics, chemistry, biology, and finance are modeled by stochastic differential equations, in which a deterministic drift is perturbed by noise \cite{KloedenPlaten1992}. Individual trajectories of such a system are random, and what is predicted and compared with data are statistics, such as expected values of observables, their relaxation toward the invariant law, and the time scales of this relaxation \cite{talay1990invariant,mattingly2002ergodicity}. The evolution of expected values is linear even when the dynamics is nonlinear. It is described by the Koopman semigroup of the process, whose generator is the Kolmogorov backward operator, a second-order differential operator built from the drift and diffusion coefficients. Approximating this semigroup on a chosen set of observables is therefore a way to predict statistics, to reduce models, and to extract time scales of a stochastic system from its coefficients or from data.

The Koopman operator \cite{Koopman1931} describes a dynamical system through the linear evolution of observables, and over the past two decades this viewpoint has grown into a family of numerical methods \cite{Mezic2005,Kutz2016}. Dynamic mode decomposition (DMD) \cite{Schmid2010,Tu2014} and extended DMD (EDMD) \cite{Williams2015} approximate the Koopman operator by a matrix. One chooses a finite dictionary of observables, evaluates it on pairs of snapshots separated by a fixed time lag, and regresses the later evaluations on the earlier ones. The matrix so obtained is a Galerkin approximation of the operator on the span of the dictionary, and it converges to the operator as the number of snapshots and then the dictionary grow \cite{KordaMezic2018}.

For stochastic systems the operator maps an observable to its conditional expectation after a time lag. This stochastic Koopman operator goes back to Mezi\'c \cite{Mezic2005}. \v{C}rnjari\'c-\v{Z}ic, Ma\'ce\v{s}i\'c, and Mezi\'c \cite{vcrnjaric2020koopman} characterized its spectrum and eigenfunctions for random dynamical systems and proposed DMD algorithms for it, and Wanner and Mezi\'c \cite{wanner2022robust} analyzed the bias that noise introduces into DMD-type approximations and proposed variants that avoid it. Snapshot-pair methods of this kind regress observables at a fixed lag. When the lag is short, the increment of an observable is dominated by its martingale part, and the variance of the differenced estimate of the generator grows like the inverse lag. When the drift and diffusion coefficients are known, the generator can instead be evaluated on the dictionary in closed form. Generator EDMD of Klus et al.~\cite{KlusEtAl2020} approximates the generator by a Galerkin matrix in this way and avoids the differencing. Finite-data error bounds for both routes are given in \cite{nuske2023finite,llamazares2024data,philipp2026variance}. For deterministic maps on the torus, EDMD with a Fourier dictionary and samples on an equispaced grid coincides with trigonometric interpolation, which gives error bounds of optimal order and a matrix-free implementation by fast Fourier transforms \cite{bartel2026extendeddynamicmodedecomposition}.

Stochastic DMD (SDMD) \cite{SDMD2025} takes the generator route and returns a discrete-time approximation of the Koopman semigroup. From samples of the invariant law it computes two matrices, the Gram matrix of the dictionary and its counterpart for the generator, and advances observables by a first-order stochastic Taylor step. We call this scheme SDMD-1. The time step is chosen after the matrices are computed, so the same matrices serve every lag, and no differencing of short-lag observations enters the estimate. In return the scheme is only first-order accurate in the time step, so that the number of steps needed for a given accuracy grows like the inverse of that accuracy. The second-order stochastic Taylor step of Appendix~A in Xu et al.~\cite{SDMD2025}, which we call SDMD-2, is the natural way to halve this exponent, but it needs the second action of the generator on every dictionary function at every sample. This requires four derivatives of the dictionary functions instead of two, the evaluation of a fourth-order differential operator in place of a second-order one, and a third matrix to compute and store.

This paper asks whether second-order accuracy can be obtained without the third matrix, and compares SDMD-1 and SDMD-2 with a third scheme, SDMD-R, obtained by Richardson extrapolation. One full step of SDMD-1 and two half steps are combined so that the leading error cancels. SDMD-R uses the matrices of SDMD-1 and nothing else. For the coefficient equations the extrapolated step coincides with a classical second-order Runge--Kutta step, so SDMD-R is not a new time integrator. As the time step tends to zero, the three schemes converge to the semigroup generated by the Galerkin generator matrix, which can also be computed from the same two matrices as a matrix exponential and serves here as the reference for the time stepping. What is new is the comparison of the three schemes as approximations of the Koopman semigroup. The purpose of the paper is to decide which second-order scheme should be used, and the answer is SDMD-R. The contribution is twofold.
\begin{enumerate}
\item On any fixed dictionary, SDMD-R is second-order accurate in the time step and needs nothing beyond SDMD-1, that is, the same two matrices, the same samples, and dictionary functions with two derivatives instead of four. The two second-order schemes converge to the same limit as the time step tends to zero, the semigroup of the Galerkin generator matrix, SDMD-R at second order and SDMD-2 in general only at first order. The difference lies in the approach to this common limit and not in the error against the Koopman semigroup, so that SDMD-R has the accuracy of SDMD-2 at the cost of SDMD-1 (\cref{thm:main}).
\item For reversible periodic diffusions with trigonometric dictionaries, the three schemes are contractions, so that no time step increases the error, as long as the time step is small enough for the size of the dictionary, and for every smooth observable their error against the Koopman semigroup is bounded explicitly and simultaneously in the size of the dictionary, the number of samples, and the time step, with the same orders for SDMD-R and SDMD-2 (\cref{thm:fourier}).
\end{enumerate}

\Cref{sec:setting} recalls the Koopman semigroup and its Galerkin approximation, reviews SDMD-1 and the Taylor scheme SDMD-2, and derives SDMD-R by Richardson extrapolation, with a comparison of cost. \Cref{sec:fixedN} contains the fixed-dictionary analysis and the finite-data bounds. \Cref{sec:fourier} treats growing trigonometric dictionaries. \Cref{sec:numerics} reports the numerical experiments, and \cref{sec:conclusion} draws the conclusions.

\section{Koopman semigroups and the SDMD schemes}\label{sec:setting}

\subsection{Koopman semigroup, generator, and Galerkin approximation}\label{subsec:background}
We consider the autonomous It\^o stochastic differential equation (SDE)
\begin{equation}\label{eq:sde}
  dX_t = b(X_t)\,dt + \sigma(X_t)\,dW_t,\qquad X_0=x\in\R^d,
\end{equation}
with drift $b\colon\R^d\to\R^d$, diffusion coefficient $\sigma\colon\R^d\to\R^{d\times d'}$, and a $d'$-dimensional Wiener process $W_t$. Periodic models are formulated directly on the torus.

\begin{assumption}\label[assumption]{ass:sde}
The SDE \eqref{eq:sde} has a unique nonexplosive Markov solution with invariant probability measure $\rho$, and the operators $(\cK^tf)(x)\coloneqq\E[f(X_t^x)]$, $t\ge0$, form a $C_0$-semigroup on the complex Hilbert space $\cH\coloneqq L^2(\rho;\C)$, $\ip{f}{g}\coloneqq \int\overline f g\,d\rho$, with generator $\cA$ and domain $\dom(\cA)$.
\end{assumption}

Jensen's inequality and invariance of $\rho$ give contractivity, $\|\cK^tf\|_{\cH}\le\|f\|_{\cH}$, so $\cA$ is dissipative, $\operatorname{Re}\ip{f}{\cA f}\le0$ for $f\in\dom(\cA)$. When the differential representation is valid, $\cA$ is the Kolmogorov backward operator $\cA f=b\cdot\nabla f+\tfrac12\operatorname{tr}(\sigma\sigma^\top\nabla^2f)$. For $f\in\dom(\cA^{p+1})$, repeated integration of the semigroup equation gives the semigroup form of the weak It\^o--Taylor expansion \cite[Chaps.~5 and~14]{KloedenPlaten1992},
\begin{equation}\label{eq:dynkin}
\begin{gathered}
  \cK^{\dt}f=\sum_{k=0}^{p}\frac{\dt^k}{k!}\cA^kf+\int_0^{\dt}\frac{(\dt-\tau)^p}{p!}\,\cK^\tau\cA^{p+1}f\,d\tau,\\
  \Big\|\cK^{\dt}f-\sum_{k=0}^{p}\frac{\dt^k}{k!}\cA^kf\Big\|_{\cH}\le\frac{\dt^{p+1}}{(p+1)!}\|\cA^{p+1}f\|_{\cH},
\end{gathered}
\end{equation}
where the bound uses contractivity of the semigroup. The SDMD schemes below are finite-dimensional versions of the truncations of \eqref{eq:dynkin}. \Cref{ass:sde} is in force throughout. The operator-theoretic arguments use only the $C_0$ contraction semigroup on $\cH$; further sampling and regularity conditions are stated where they enter.

The approximation space is spanned by a dictionary. Let $\{\psi_1,\dots,\psi_N\}\subset\dom(\cA)$ be linearly independent and let
\[
  \cV_N\coloneqq \Span\{\psi_1,\dots,\psi_N\},\qquad \Phi_N\colon\C^N\to\cV_N,\qquad \Phi_N\mathbf v\coloneqq \sum_{j=1}^N(\mathbf v)_j\psi_j,
\]
where the synthesis map $\Phi_N$ identifies coefficient vectors with elements of $\cV_N$, and let $\Proj_N\colon\cH\to\cV_N$ be the $L^2(\rho)$-orthogonal projection. The population Gram matrix and its counterpart for the generator are
\[
  \mathbf G_{ij}\coloneqq \ip{\psi_i}{\psi_j},\qquad \mathbf H_{ij}\coloneqq \ip{\psi_i}{\cA\psi_j},
\]
$\mathbf G$ is invertible because the $\psi_j$ are linearly independent, and the Galerkin generator matrix and the Galerkin generator are
\[
  \mathbf A_N\coloneqq \mathbf G^{-1}\mathbf H,\qquad \cA_N\coloneqq \Phi_N\mathbf A_N\Phi_N^{-1}=\Proj_N\cA|_{\cV_N},
\]
so that $\ip{u}{\cA_Nv}=\ip{u}{\cA v}$ for $u,v\in\cV_N$. The matrix $\mathbf A_N$ is the large-sample limit of the generator EDMD matrix of Klus et al.~\cite{KlusEtAl2020}. If $\cA$ is self-adjoint, so is $\cA_N$ on $\cV_N$. We also use the projected Koopman step $\cK_N^{\dt}\coloneqq \Proj_N\cK^{\dt}|_{\cV_N}$, the large-sample limit of the snapshot EDMD step at lag $\dt$. Operators are written in calligraphic type and matrices and coefficient vectors in boldface, so that $\cA_N$ is an operator on $\cV_N$, $\mathbf A_N$ is its matrix, and $v=\Phi_N\mathbf v$; $I$ and $\mathbf I$ denote the identity operator and the identity matrix.

A norm without a subscript is the norm of $\cH$ or the operator norm that it induces. On $\cV_N$ it is inherited from $L^2(\rho)$, that is, $\|\Phi_N\mathbf v\|_{\cH}^2=\mathbf v^*\mathbf G\mathbf v$, and the operator norms of $\cA_N$, $e^{t\cA_N}$, and of the one-step maps are $\cV_N\to\cV_N$ norms. For a coefficient matrix $\mathbf M$,
\begin{equation}\label{eq:transport}
  \|\Phi_N\mathbf M\Phi_N^{-1}\|_{\cV_N\to\cV_N}=\|\mathbf G^{1/2}\mathbf M\mathbf G^{-1/2}\|_2\le\kappa_N\|\mathbf M\|_2,
\end{equation}
where $\|\cdot\|_2$ is the spectral norm of a matrix (and the Euclidean norm of a vector) and $\kappa_N\coloneqq \sqrt{\operatorname{cond}_2(\mathbf G)}$. Dissipativity is inherited, since for $v=\Phi_N\mathbf v$ one has $\operatorname{Re}(\mathbf v^*\mathbf G\mathbf A_N\mathbf v)=\operatorname{Re}\ip{v}{\cA v}\le0$, so that
\begin{equation}\label{eq:AN_contraction}
  \|e^{t\cA_N}\|\le1,\qquad t\ge0 .
\end{equation}
The Galerkin semigroup $e^{t\cA_N}$ is the evolution within $\cV_N$ generated by $\cA_N$; all SDMD schemes below approximate it, and \cref{sec:fourier} relates it to $\cK^t$.

\subsection{The first-order scheme SDMD-1}\label{subsec:sdmd1}
SDMD \cite{SDMD2025} approximates $\cK^{\dt}$ on $\cV_N$ by the projection of the first-order truncation of \eqref{eq:dynkin},
\begin{equation}\label{eq:sdmd1}
  \cQ^{[1]}_{N,\dt}\coloneqq \Proj_N(I+\dt\cA)|_{\cV_N}=I+\dt \cA_N ,
  \qquad \cQ^{[1]}_{\dt}(\cM)\coloneqq I+\dt \cM ,
\end{equation}
where the second formula defines the one-step map for an arbitrary operator $\cM$ on $\cV_N$. In practice $\cA_N$ is replaced by an empirical Galerkin generator. Given samples $x_1,\dots,x_m$ of $\rho$, let the rows of $\boldsymbol\Psi\in\C^{m\times N}$ contain $\psi(x_k)^\top$, where $\psi(x)\coloneqq (\psi_1(x),\dots,\psi_N(x))^\top$, and let the array $\cA\boldsymbol\Psi$ contain the evaluations $(\cA\psi_j)(x_k)$, which are available in closed form when $b$ and $\sigma$ are known. The empirical matrices, the empirical Galerkin generator with a Tikhonov penalty $\boldsymbol\Lambda_N\succeq0$ and $\gamma\ge0$, and the empirical step are
\begin{equation}\label{eq:empirical_matrices}
  \widehat{\mathbf G}_m\coloneqq \frac1m\boldsymbol\Psi^*\boldsymbol\Psi,\qquad
  \widehat{\mathbf H}_m\coloneqq \frac1m\boldsymbol\Psi^*(\cA\boldsymbol\Psi),\qquad
  \widehat{\mathbf A}_{N,m,\gamma}\coloneqq (\widehat{\mathbf G}_m+\gamma\boldsymbol\Lambda_N)^{-1}\widehat{\mathbf H}_m,
\end{equation}
$\widehat{\cA}_{N,m,\gamma}\coloneqq \Phi_N\widehat{\mathbf A}_{N,m,\gamma}\Phi_N^{-1}$, and $\widehat{\cQ}^{[1]}_{N,\dt,m,\gamma}\coloneqq \cQ^{[1]}_{\dt}(\widehat{\cA}_{N,m,\gamma})$. In coefficients, an observable $v=\Phi_N\mathbf v$ is advanced to the time $T=n\dt$ by $\mathbf v\mapsto\mathbf Q^n\mathbf v$ with $\mathbf Q\coloneqq \mathbf I+\dt\widehat{\mathbf A}_{N,m,\gamma}$, either by $n$ matrix--vector products or by one matrix power. The time step enters only after the matrices are computed, so the same matrices serve every lag. The scheme is first-order accurate in $\dt$ (\cref{lem:onestep} below). Positivity and the Markov property are not imposed on this or on the following polynomial propagators; stability is treated separately. The relation to the statements of Xu et al.~\cite{SDMD2025} is discussed in \cref{app:relation}.

\subsection{The second-order Taylor scheme SDMD-2}\label{subsec:sdmd2}
The same route gives a second-order scheme, the one of Appendix~A in Xu et al.~\cite{SDMD2025}. If $\cV_N\subset\dom(\cA^2)$, the projection of the second-order truncation of \eqref{eq:dynkin} is
\begin{equation}\label{eq:sdmd2}
  \cQ^{[2]}_{N,\dt}\coloneqq \Proj_N\Big(I+\dt\cA+\frac{\dt^2}{2}\cA^2\Big)\Big|_{\cV_N}=I+\dt \cA_N+\frac{\dt^2}{2}\cB_N,
\end{equation}
where $\cB_N\coloneqq \Proj_N\cA^2|_{\cV_N}$, with the one-step map $\cQ^{[2]}_{\dt}(\cM,\cN)\coloneqq I+\dt \cM+\tfrac12\dt^2\cN$ for operators $\cM$, $\cN$ on $\cV_N$. The new operator has the matrix representation
\[
  \mathbf J_{ij}\coloneqq \ip{\psi_i}{\cA^2\psi_j},\qquad \mathbf B_N\coloneqq \mathbf G^{-1}\mathbf J,\qquad \cB_N=\Phi_N\mathbf B_N\Phi_N^{-1},
\]
and its empirical version requires the array $\cA^2\boldsymbol\Psi$ of the evaluations $(\cA^2\psi_j)(x_k)$,
\[
\begin{gathered}
  \widehat{\mathbf J}_m\coloneqq \frac1m\boldsymbol\Psi^*(\cA^2\boldsymbol\Psi),\qquad
  \widehat{\mathbf B}_{N,m,\gamma}\coloneqq (\widehat{\mathbf G}_m+\gamma\boldsymbol\Lambda_N)^{-1}\widehat{\mathbf J}_m,\\
  \widehat{\cB}_{N,m,\gamma}\coloneqq \Phi_N\widehat{\mathbf B}_{N,m,\gamma}\Phi_N^{-1},
\end{gathered}
\]
with the empirical step $\widehat{\cQ}^{[2]}_{N,\dt,m,\gamma}\coloneqq \cQ^{[2]}_{\dt}(\widehat{\cA}_{N,m,\gamma},\widehat{\cB}_{N,m,\gamma})$. The second order comes at three costs. The dictionary must lie in $\dom(\cA^2)$, which for a second-order differential operator $\cA$ means four derivatives. The array $\cA^2\boldsymbol\Psi$ requires $mN$ evaluations of the fourth-order operator $\cA^2$, which in $d$ dimensions involve derivatives of $\psi_j$ up to order four and of $b$ and $\sigma\sigma^\top$ up to order two ($O(d^4)$ terms against $O(d^2)$ for $\cA\psi_j$). The third matrix $\widehat{\mathbf J}_m$ costs a further $mN^2$ multiply--adds and carries its own sampling error. For self-adjoint generators an energy form avoids the evaluations of $\cA^2\psi_j$ and the two additional derivatives (\cref{app:energy}), but not the third matrix; it is tested on a spline dictionary in \cref{subsec:num_spline}.

\subsection{Richardson extrapolation and the scheme SDMD-R}\label{subsec:sdmdR}
Richardson extrapolation \cite{Richardson1911}, \cite[Sec.~II.9]{HairerNorsettWanner1993} raises the order of an approximation whose error has an expansion in the step size. If the local error of a one-step map $\cQ_{\dt}$ has a leading term of order $\dt^2$, then two half steps have one half of this leading term, and the combination $2\cQ_{\dt/2}^2-\cQ_{\dt}$ cancels it. Applied to the SDMD-1 map this gives
\begin{equation}\label{eq:sdmdR}
  \QR_{\dt}(\cM)\coloneqq 2\Big(I+\frac{\dt}{2}\cM\Big)^2-(I+\dt \cM)=I+\dt \cM+\frac{\dt^2}{2}\cM^2,
\end{equation}
with the population step $\QR_{N,\dt}\coloneqq \QR_{\dt}(\cA_N)$ and the empirical step $\widehat{\cQ}^{\mathrm{R}}_{N,\dt,m,\gamma}\coloneqq \QR_{\dt}(\widehat{\cA}_{N,m,\gamma})$. We call this scheme SDMD-R. It uses the matrices $\widehat{\mathbf G}_m$ and $\widehat{\mathbf H}_m$ of SDMD-1 and nothing else, requires $\cV_N\subset\dom(\cA)$ only, forms $\widehat{\mathbf A}_{N,m,\gamma}^2$ once at cost $O(N^3)$, and then steps at the cost of SDMD-1. For the linear equation $\dot{\mathbf v}=\mathbf A_N\mathbf v$, every explicit two-stage Runge--Kutta method of order two \cite[Sec.~II.1]{HairerNorsettWanner1993} advances $\mathbf v$ by $\QR_{\dt}(\mathbf A_N)$. The maps $\cQ^{[1]}_{\dt}(\cM)$ and $\QR_{\dt}(\cM)$ are the Taylor polynomials of degrees one and two of the exponential $e^{\dt \cM}$, and $e^{\dt \cA_N}$ is the exact step of the Galerkin semigroup of \cref{subsec:background}. With the matrices of SDMD-1 this step can also be evaluated as a matrix exponential, at a cost of order $N^3$ operations for each step size up to a factor logarithmic in the norm of the matrix \cite{Higham2005}, and polynomials of higher degree in the generator matrix give higher orders from the same two matrices. We keep to degree two, the order of SDMD-2, and use the exponential as the reference for the time stepping (\cref{rem:exponential}). \Cref{tab:cost} compares the three schemes. When $m\gg N$ and $m\gg n$, the computation of the matrices dominates the total cost.

\begin{table}[tbp]
\footnotesize
\centering
\setlength{\tabcolsep}{3.5pt}
\begin{tabular}{L{0.11\linewidth}L{0.14\linewidth}L{0.16\linewidth}L{0.12\linewidth}L{0.17\linewidth}L{0.17\linewidth}}
\toprule
Scheme & Evaluations per sample & Matrices & Flops for the matrices & Regularity & Time stepping \\
\midrule
SDMD-1 & $\psi$, $\cA\psi$ & $\widehat{\mathbf G}_m$, $\widehat{\mathbf H}_m$ & $2mN^2$ & $\cV_N\subset\dom(\cA)$ & $O(N^3+nN^2)$ \\
SDMD-R & $\psi$, $\cA\psi$ & $\widehat{\mathbf G}_m$, $\widehat{\mathbf H}_m$ & $2mN^2$ & $\cV_N\subset\dom(\cA)$ & $O(N^3+nN^2)$ \\
SDMD-2 & $\psi$, $\cA\psi$, $\cA^2\psi$ & $\widehat{\mathbf G}_m$, $\widehat{\mathbf H}_m$, $\widehat{\mathbf J}_m$ & $3mN^2$ & $\cV_N\subset\dom(\cA^2)$ & $O(N^3+nN^2)$ \\
\bottomrule
\end{tabular}
\captionbelow{Requirements of the three schemes for $m$ samples, $N$ dictionary functions, and $n=T/\dt$ steps; flop counts are leading-order multiply--adds and exclude the evaluation of $\cA\psi$ and $\cA^2\psi$.}
\label{tab:cost}
\end{table}

The two second-order schemes differ in the term of order $\dt^2$. SDMD-2 projects the second generator action, SDMD-R squares the projected generator, and the difference is the \emph{second-order Galerkin residual}
\[
\begin{gathered}
  \cD_N\coloneqq \cB_N-\cA_N^2,\qquad \Theta_N\coloneqq \|\cD_N\|
  =\big\|\mathbf G^{1/2}(\mathbf B_N-\mathbf A_N^2)\mathbf G^{-1/2}\big\|_2,\\
  \QR_{N,\dt}-\cQ^{[2]}_{N,\dt}=-\frac{\dt^2}{2}\cD_N .
\end{gathered}
\]
The operator $\cE_N^\perp\coloneqq (I-\Proj_N)\cA|_{\cV_N}$ is the generator invariance defect of the dictionary, and the residual is the part of this defect that returns to $\cV_N$ after a second application of $\cA$.

\begin{proposition}\label[proposition]{prop:residual_structure}
Let $\cV_N\subset\dom(\cA^2)$. Then $\cD_N=\Proj_N\cA(I-\Proj_N)\cA|_{\cV_N}$. If $\cA$ is self-adjoint, then $\cD_N=(\cE_N^\perp)^*\cE_N^\perp\succeq0$ and $\Theta_N=\|\cE_N^\perp\|^2$; in particular $\cD_N=0$ if and only if $\cA\cV_N\subseteq\cV_N$.
\end{proposition}

\begin{proof}
For $v\in\cV_N$, $\cA v\in\dom(\cA)$ and $\Proj_N\cA v\in\cV_N\subset\dom(\cA)$, so $(I-\Proj_N)\cA v\in\dom(\cA)$ and $\cD_Nv=\Proj_N\cA^2v-\Proj_N\cA\Proj_N\cA v=\Proj_N\cA(I-\Proj_N)\cA v$. If $\cA$ is self-adjoint, then for $u,v\in\cV_N$, $\ip{u}{\cD_Nv}=\ip{\cA u}{(I-\Proj_N)\cA v}=\ip{\cE_N^\perp u}{\cE_N^\perp v}$, because $\Proj_N\cA u\perp(I-\Proj_N)\cA v$.
\end{proof}

The rest of the paper quantifies the effect of the invariance defect, which enters the time-stepping error of SDMD-2 through $\cD_N$ and the sampling error of every scheme through $\cE_N^\perp$.

\section{Error analysis at a fixed dictionary}\label{sec:fixedN}

\subsection{One-step defects and stability}\label{subsec:onestep}
Formula \eqref{eq:dynkin} holds for the semigroup $e^{t\cA_N}$ on $\cV_N$, for every $v\in\cV_N$ and every $p$, with the bound $\dt^{p+1}\|\cA_N\|^{p+1}\|v\|/(p+1)!$ by \eqref{eq:AN_contraction}.

\begin{lemma}\label[lemma]{lem:onestep}
Let $\cV_N\subset\dom(\cA)$ and $\dt>0$. Then
\begin{equation}\label{eq:onestep_1R}
  \|e^{\dt \cA_N}-\cQ^{[1]}_{N,\dt}\|\le\tfrac12\|\cA_N\|^2\dt^2,\qquad
  \|e^{\dt \cA_N}-\QR_{N,\dt}\|\le\tfrac16\|\cA_N\|^3\dt^3 .
\end{equation}
If $\cV_N\subset\dom(\cA^2)$, then
\begin{equation}\label{eq:onestep_2}
\begin{gathered}
  \|e^{\dt \cA_N}-\cQ^{[2]}_{N,\dt}\|\le\tfrac12\Theta_N\dt^2+\tfrac16\|\cA_N\|^3\dt^3,\\
  \lim_{\dt\downarrow0}\dt^{-2}\|e^{\dt \cA_N}-\cQ^{[2]}_{N,\dt}\|=\tfrac12\Theta_N .
\end{gathered}
\end{equation}
\end{lemma}

\begin{proof}
Formula \eqref{eq:dynkin} for $e^{t\cA_N}$ with $p=1$ and $p=2$ gives \eqref{eq:onestep_1R}, because $\cQ^{[1]}_{N,\dt}$ and $\QR_{N,\dt}$ are the Taylor polynomials of $e^{\dt \cA_N}$ of degrees one and two. Adding and subtracting $\tfrac12\dt^2\cA_N^2$ and using $\|\cB_N-\cA_N^2\|=\Theta_N$ gives the bound in \eqref{eq:onestep_2}; the limit follows from $e^{\dt \cA_N}-\cQ^{[2]}_{N,\dt}=\tfrac12\dt^2(\cA_N^2-\cB_N)+O(\dt^3)$.
\end{proof}

Fix $T>0$ and $\dt_0>0$. A family of one-step maps $\cQ_{\dt}$ on $\cV_N$ is called stable on $[0,T]$ if $\sup_{0\le k\le T/\dt}\|\cQ_{\dt}^k\|\le M_T$ for all $0<\dt\le\dt_0$ with $T/\dt\in\mathbb N$. At fixed $N$ each of the maps \eqref{eq:sdmd1}, \eqref{eq:sdmd2}, and \eqref{eq:sdmdR} is stable, since $\|\cQ_{N,\dt}\|\le1+\dt\|\cA_N\|+\tfrac12\dt^2C_N\le\exp[\dt(\|\cA_N\|+\tfrac12\dt_0C_N)]$ with $C_N=0$, $\|\cB_N\|$, or $\|\cA_N\|^2$, so that $M_T=\exp[T(\|\cA_N\|+\tfrac12\dt_0C_N)]$. For bounded operators $\cS$ and $\cQ$,
\begin{equation}\label{eq:telescoping}
  \cS^n-\cQ^n=\sum_{j=0}^{n-1}\cS^{n-1-j}(\cS-\cQ)\cQ^j,
  \qquad
  \|\cS^n-\cQ^n\|\le nM_SM_Q\|\cS-\cQ\|
\end{equation}
when $\sup_{k\le n}\|\cS^k\|\le M_S$ and $\sup_{k\le n}\|\cQ^k\|\le M_Q$. If $\cQ_{\dt}$ is stable with constant $M_T$ and $\|\widehat{\cQ}_{\dt}-\cQ_{\dt}\|\le\alpha\dt$, then \eqref{eq:telescoping} and the discrete Gr\"onwall inequality give
\begin{equation}\label{eq:stability_transfer}
  \sup_{0\le k\le T/\dt}\|\widehat{\cQ}_{\dt}^k\|\le M_T\exp(M_T\alpha T),
\end{equation}
which transfers stability from the population map to a perturbed one.

\begin{lemma}\label[lemma]{lem:insertion}
For operators $\cA_1,\cA_2,\cB_1,\cB_2$ on $\cV_N$,
\[
\begin{aligned}
  \|\cQ^{[1]}_{\dt}(\cA_1)-\cQ^{[1]}_{\dt}(\cA_2)\|&\le\dt\|\cA_1-\cA_2\|,\\
  \|\cQ^{[2]}_{\dt}(\cA_1,\cB_1)-\cQ^{[2]}_{\dt}(\cA_2,\cB_2)\|&\le\dt\|\cA_1-\cA_2\|+\tfrac12\dt^2\|\cB_1-\cB_2\|,\\
  \|\QR_{\dt}(\cA_1)-\QR_{\dt}(\cA_2)\|&\le\dt\big(1+\tfrac12\dt(\|\cA_1\|+\|\cA_2\|)\big)\|\cA_1-\cA_2\| .
\end{aligned}
\]
\end{lemma}
\begin{proof}
Subtract and use $\cA_1^2-\cA_2^2=\cA_1(\cA_1-\cA_2)+(\cA_1-\cA_2)\cA_2$.
\end{proof}

\subsection{Fixed-horizon error}\label{subsec:main}
The following theorem is the first main result. It is stated for arbitrary perturbations of the Galerkin matrices, so that exact matrices ($\xi_A=\xi_B=0$), that is, the population matrices of \cref{sec:setting}, and empirical matrices (\cref{prop:perturbation} below) are covered by one statement.

\begin{theorem}[Fixed-horizon error]\label[theorem]{thm:main}
Let $\cV_N\subset\dom(\cA)$, $T>0$, and $0<\dt\le\dt_0$ with $n\coloneqq T/\dt\in\mathbb N$. Let $\widehat{\cA}$ be an operator on $\cV_N$ with $\|\widehat{\cA}-\cA_N\|\le\xi_A\le\overline\xi_A$ and, for SDMD-2, let $\cV_N\subset\dom(\cA^2)$ and $\|\widehat{\cB}-\cB_N\|\le\xi_B\le\overline\xi_B$. Write $\widehat{\cQ}^{[1]}\coloneqq \cQ^{[1]}_{\dt}(\widehat{\cA})$, $\widehat{\cQ}^{\mathrm R}\coloneqq \QR_{\dt}(\widehat{\cA})$, $\widehat{\cQ}^{[2]}\coloneqq \cQ^{[2]}_{\dt}(\widehat{\cA},\widehat{\cB})$, and let
\[
  d_N(T,f)\coloneqq \|e^{T\cA_N}\Proj_Nf-\cK^Tf\|_{\cH}
\]
be the Galerkin error. Then, for every $f\in\cH$,
\begin{equation}\label{eq:main}
\begin{aligned}
  \|(\widehat{\cQ}^{[1]})^n\Proj_Nf-\cK^Tf\|_{\cH}&\le d_N(T,f)+C(\dt+\xi_A)\|f\|_{\cH},\\
  \|(\widehat{\cQ}^{\mathrm R})^n\Proj_Nf-\cK^Tf\|_{\cH}&\le d_N(T,f)+C(\dt^2+\xi_A)\|f\|_{\cH},\\
  \|(\widehat{\cQ}^{[2]})^n\Proj_Nf-\cK^Tf\|_{\cH}&\le d_N(T,f)+C(\Theta_N\dt+\dt^2+\xi_A+\dt\,\xi_B)\|f\|_{\cH},
\end{aligned}
\end{equation}
where $C$ depends only on $T$, $\dt_0$, $\|\cA_N\|$, $\|\cB_N\|$, $\overline\xi_A$, and $\overline\xi_B$. Moreover, the term $\Theta_N\dt$ in \eqref{eq:main} is attained in the sense that, with exact matrices,
\begin{equation}\label{eq:leading_term}
  \Big\|(\cQ^{[2]}_{N,\dt})^n-e^{T\cA_N}-\dt\,\cL_{N,T}\Big\|\le C\dt^2,
  \quad
  \cL_{N,T}\coloneqq \frac12\int_0^Te^{(T-\tau)\cA_N}\cD_Ne^{\tau \cA_N}\,d\tau ,
\end{equation}
and, if $\cA$ is self-adjoint and $\{(\lambda_j,w_j)\}_{j=1}^N$ is an orthonormal eigenbasis of $\cA_N$,
\begin{equation}\label{eq:lower_bound}
  \|\cL_{N,T}\|\ge\frac{T}{2}\max_{1\le j\le N}e^{T\lambda_j}\ip{w_j}{\cD_Nw_j},
\end{equation}
which is positive unless $\cA\cV_N\subseteq\cV_N$; in that case $\cD_N=0$ and SDMD-2 coincides with SDMD-R.
\end{theorem}

\begin{proof}
For each scheme, with $\cQ$ the population map and $\widehat{\cQ}$ the perturbed one,
\[
  \widehat{\cQ}^n\Proj_Nf-\cK^Tf
  =\big[\widehat{\cQ}^n-\cQ^n\big]\Proj_Nf+\big[\cQ^n-e^{T\cA_N}\big]\Proj_Nf+\big[e^{T\cA_N}\Proj_Nf-\cK^Tf\big],
\]
and the last term has norm $d_N(T,f)$. For the second term, $e^{\dt \cA_N}$ is a contraction by \eqref{eq:AN_contraction} and $\cQ$ is stable with constant $M_T$, so \eqref{eq:telescoping} with $\cS=e^{\dt \cA_N}$ and \cref{lem:onestep} give $\|\cQ^n-e^{T\cA_N}\|\le nM_T\|\cQ-e^{\dt \cA_N}\|$, which is at most $\tfrac12TM_T\|\cA_N\|^2\dt$ for SDMD-1, $\tfrac16TM_T\|\cA_N\|^3\dt^2$ for SDMD-R, and $TM_T(\tfrac12\Theta_N\dt+\tfrac16\|\cA_N\|^3\dt^2)$ for SDMD-2. For the first term, \cref{lem:insertion} gives $\|\widehat{\cQ}-\cQ\|\le\alpha\dt$ with $\alpha=\xi_A$ for SDMD-1, $\alpha=\xi_A+\tfrac12\dt\,\xi_B$ for SDMD-2, and $\alpha=(1+\dt_0(\|\cA_N\|+\tfrac12\overline\xi_A))\xi_A$ for SDMD-R, because $\|\widehat{\cA}\|\le\|\cA_N\|+\overline\xi_A$. Then \eqref{eq:stability_transfer} bounds the powers of $\widehat{\cQ}$ by $\widehat M_T\coloneqq M_T\exp(M_T\overline\alpha T)$ with $\overline\alpha$ the value of $\alpha$ at $\xi_A=\overline\xi_A$, $\xi_B=\overline\xi_B$, and $\dt=\dt_0$, and \eqref{eq:telescoping} gives $\|\widehat{\cQ}^n-\cQ^n\|\le nM_T\widehat M_T\alpha\dt=TM_T\widehat M_T\alpha$. Collecting the three terms proves \eqref{eq:main}.

For \eqref{eq:leading_term}, let $\cS\coloneqq e^{\dt \cA_N}$ and $\cQ\coloneqq \cQ^{[2]}_{N,\dt}$. By \cref{lem:onestep}, $\cS=\QR_{N,\dt}+\dt^3\cZ$ with $\|\cZ\|\le\tfrac16\|\cA_N\|^3$, so that $\cQ=\cS+\tfrac12\dt^2\cD_N-\dt^3\cZ$. Since $\|\cS\|\le1$ and $\|\cQ\|\le1+\tfrac12\dt^2\Theta_N+\dt^3\|\cZ\|$, the powers $\cQ^k$, $k\le n$, are bounded by a constant $M$ independent of $\dt$, and $\|\cQ^k-\cS^k\|\le kM\|\cQ-\cS\|\le C\dt$ by \eqref{eq:telescoping}. The identity \eqref{eq:telescoping} with the roles of $\cS$ and $\cQ$ exchanged gives
\[
\begin{aligned}
  \cQ^n-\cS^n
  &=\frac{\dt^2}{2}\sum_{j=0}^{n-1}\cS^{n-1-j}\cD_N\cS^j
   +\frac{\dt^2}{2}\sum_{j=0}^{n-1}\big(\cQ^{n-1-j}-\cS^{n-1-j}\big)\cD_N\cS^j\\
  &\quad-\dt^3\sum_{j=0}^{n-1}\cQ^{n-1-j}\cZ \cS^j .
\end{aligned}
\]
The norms of the second and the third term are at most
\[
  \tfrac12\dt^2\,n\,C\dt\,\Theta_N\le C\dt^2
  \quad\text{and}\quad
  \dt^3nM\|\cZ\|\le C\dt^2 .
\]
In the first term, $\cS^{n-1-j}=e^{-\dt \cA_N}e^{(T-j\dt)\cA_N}$, so this term is $\dt\,e^{-\dt \cA_N}$ times a Riemann sum with step $\dt$ for the integral that defines $\cL_{N,T}$, whose integrand is Lipschitz continuous. It therefore equals
\[
  \dt\,e^{-\dt \cA_N}\big(\cL_{N,T}+O(\dt)\big)=\dt\,\cL_{N,T}+O(\dt^2),
\]
because $\|e^{-\dt \cA_N}-I\|\le C\dt$. This proves \eqref{eq:leading_term}. For \eqref{eq:lower_bound}, \cref{prop:residual_structure} gives $\cD_N\succeq0$ and, since $\cA_Nw_j=\lambda_jw_j$ with $\lambda_j$ real,
\[
  \ip{w_j}{\cL_{N,T}w_j}=\frac12\int_0^Te^{(T-\tau)\lambda_j}e^{\tau\lambda_j}\ip{w_j}{\cD_Nw_j}\,d\tau=\frac{T}{2}e^{T\lambda_j}\ip{w_j}{\cD_Nw_j}.
\]
The operator norm dominates every quadratic form, and $\sum_j\ip{w_j}{\cD_Nw_j}=\operatorname{tr}\cD_N>0$ unless $\cD_N=0$, which by \cref{prop:residual_structure} means $\cA\cV_N\subseteq\cV_N$.
\end{proof}

\begin{remark}\label[remark]{rem:reading}
The Galerkin error $d_N(T,f)$ is the error of the Galerkin semigroup $e^{T\cA_N}$, to which all three schemes converge as $\dt\to0$. It does not depend on the scheme and is the subject of \cref{sec:fourier}. The scheme-dependent part consists of the time-stepping term and the data term, and SDMD-R reaches second order under the hypothesis $\cV_N\subset\dom(\cA)$ of SDMD-1 and with the same perturbation $\xi_A$. The term $\Theta_N\dt$ separates the two second-order schemes in their distance from the Galerkin semigroup. It is not an error against the Koopman semigroup. With exact matrices the error of either scheme against $\cK^Tf$ tends to $d_N(T,f)$ as $\dt\to0$, that of SDMD-R at second order and that of SDMD-2 at first order with the leading operator $\dt\,\cL_{N,T}$, and for $\dt>0$ the error of SDMD-2 can be smaller or larger than that of SDMD-R (\cref{sec:numerics}). By \eqref{eq:lower_bound}, the invariance defect matters at the horizon $T$ when it sits on slowly decaying eigenvectors of $\cA_N$ and is damped by $e^{T\lambda_j}$ when it sits on rapidly decaying ones. The Ornstein--Uhlenbeck (OU) example of \cref{subsec:num_ou} is of the first kind, the Fourier example of \cref{subsec:num_fourier} of the second. The same leading term governs the iterated projected Koopman step, since for $\cV_N\subset\dom(\cA^3)$, \eqref{eq:dynkin} with $p=2$ gives $\cK_N^{\dt}=e^{\dt \cA_N}+\tfrac12\dt^2\cD_N+O(\dt^3)$ and the proof of \eqref{eq:leading_term} applies verbatim to $(\cK_N^{\dt})^n$. To first order in $\dt$, SDMD-2 thus follows the iterated projected Koopman step, which for $n=1$ maps each $v\in\cV_N$ to the element of $\cV_N$ closest to $\cK^Tv$, and SDMD-R follows the Galerkin semigroup. In the data term, SDMD-2 carries $\dt\,\xi_B$ in addition to $\xi_A$, the contribution of the matrix that SDMD-R does not form. \Cref{prop:perturbation} shows that $\xi_B$ has the same sampling rate as $\xi_A$, and \cref{sec:fourier} quantifies both in the cutoff.
\end{remark}

\begin{remark}\label[remark]{rem:exponential}
The limit $e^{T\cA_N}$ of the three schemes can be computed from the matrices of SDMD-1 as a matrix exponential. For $\widehat{\cA}$ as in \cref{thm:main}, dissipativity of $\cA_N$ gives $\|e^{t\widehat{\cA}}\|\le e^{t\xi_A}$ for $t\ge0$, and Duhamel's formula $e^{T\widehat{\cA}}-e^{T\cA_N}=\int_0^Te^{(T-\tau)\cA_N}(\widehat{\cA}-\cA_N)e^{\tau\widehat{\cA}}\,d\tau$ with \eqref{eq:AN_contraction} gives
\begin{equation}\label{eq:exponential}
  \|e^{T\widehat{\cA}}\Proj_Nf-\cK^Tf\|_{\cH}\le d_N(T,f)+T\xi_Ae^{T\xi_A}\|f\|_{\cH} .
\end{equation}
This bound has no time-stepping term, and its constant does not depend on $\|\cA_N\|$, whereas the constant $C$ of \eqref{eq:main} depends on $\|\cA_N\|$, exponentially through the stability constant $M_T$ unless the one-step maps are contractions, as in \cref{sec:fourier}. The exponential requires a function of the full matrix, and the polynomial schemes only products with the generator matrix. By \eqref{eq:main}, the bound for SDMD-R is of the order of the bound for the exponential as soon as $\dt^2$ is below the Galerkin and data terms, and the bound for SDMD-1 as soon as $\dt$ is. Both cases, a time-stepping error above and below the Galerkin error, occur in the experiments of \cref{sec:numerics}.
\end{remark}

\subsection{Empirical Galerkin matrices}\label{subsec:data}
The perturbations $\xi_A$ and $\xi_B$ of \cref{thm:main} are now quantified for the empirical matrices \eqref{eq:empirical_matrices} under the following hypotheses.
\begin{enumerate}
\item[(H1)] The samples $x_1,\dots,x_m$ are independent and identically distributed with law $\rho$; the dictionary is deterministic or independent of the samples, and all statements are then conditional on the dictionary.
\item[(H2)] There are constants $F_0,F_1,g_*>0$ with $\|\psi(x)\|_2\le F_0$ and $\|\cA\psi(x)\|_2\le F_1$ for $\rho$-almost every $x$, and $\mathbf G\succeq g_*\mathbf I$. For SDMD-2, in addition $\cV_N\subset\dom(\cA^2)$ and $\|\cA^2\psi(x)\|_2\le F_2$.
\end{enumerate}
Independence in (H1) enters through the matrix Bernstein inequality. Samples from a trajectory or from the invariant law of a numerical integrator are discussed in \cref{rem:data_variants}.

The Galerkin generator leaves a residual on each dictionary function. Define
\[
  r_j\coloneqq \cA\psi_j-\cA_N\psi_j=\cE_N^\perp\psi_j,\qquad
  r(x)\coloneqq (\cA\psi)(x)-\mathbf A_N^\top\psi(x),
\]
so that $r_j\in\cV_N^\perp$ by the Galerkin relation, and, when $\cV_N\subset\dom(\cA^2)$, $s_j\coloneqq \cA^2\psi_j-\cB_N\psi_j=(I-\Proj_N)\cA^2\psi_j$ and $s(x)\coloneqq (\cA^2\psi)(x)-\mathbf B_N^\top\psi(x)$. For the high-probability bounds, set
\[
\begin{gathered}
  \eta_{m,N,\delta}\coloneqq \sqrt{\frac{\log(2N/\delta)}{m}}+\frac{\log(2N/\delta)}{m},\\
  F_r\coloneqq \operatorname*{ess\,sup}_x\|r(x)\|_2,
  \qquad
  F_s\coloneqq \operatorname*{ess\,sup}_x\|s(x)\|_2,
\end{gathered}
\]
with $F_r\le F_1+\|\mathbf A_N\|_2F_0$ and $F_s\le F_2+\|\mathbf B_N\|_2F_0$.

\begin{lemma}\label[lemma]{lem:residual}
For arbitrary samples $x_1,\dots,x_m$,
\begin{equation}\label{eq:residual_identity}
\begin{aligned}
  \widehat{\mathbf H}_m-\widehat{\mathbf G}_m\mathbf A_N&=\mathbf R_m\coloneqq \frac1m\sum_{k=1}^m\overline{\psi(x_k)}\,r(x_k)^\top,\\
  \widehat{\mathbf A}_{N,m,\gamma}-\mathbf A_N&=(\widehat{\mathbf G}_m+\gamma\boldsymbol\Lambda_N)^{-1}\big(\mathbf R_m-\gamma\boldsymbol\Lambda_N\mathbf A_N\big)
\end{aligned}
\end{equation}
whenever $\widehat{\mathbf G}_m+\gamma\boldsymbol\Lambda_N$ is invertible, and $\E_\rho[\overline{\psi(x)}\,r(x)^\top]=0$; the same holds for $\widehat{\mathbf J}_m$, $\mathbf B_N$, and $s$. If $\cA\cV_N\subseteq\cV_N$, then $r=0$ $\rho$-almost everywhere, and for samples of $\rho$, exact generator evaluations, $\gamma=0$, and invertible $\widehat{\mathbf G}_m$, $\widehat{\mathbf A}_{N,m,0}=\mathbf A_N$ almost surely.
\end{lemma}

\begin{proof}
Since $\cA_N\psi_j=\sum_l(\mathbf A_N)_{lj}\psi_l$, the $(i,j)$ entry of $\widehat{\mathbf G}_m\mathbf A_N$ is
\[
  \frac1m\sum_{k=1}^m\overline{\psi_i(x_k)}\,(\cA_N\psi_j)(x_k),
\]
which gives the identity; subtracting $(\widehat{\mathbf G}_m+\gamma\boldsymbol\Lambda_N)\mathbf A_N$ from $\widehat{\mathbf H}_m$ gives the second formula. The expectation vanishes because $\E_\rho[\overline{\psi_i}\,r_j]=\ip{\psi_i}{r_j}$ and $r_j\perp\cV_N$. For an invariant dictionary $r_j=\cE_N^\perp\psi_j=0$ in $\cH$, so $\mathbf R_m=0$ almost surely.
\end{proof}

For $\gamma=0$ the empirical generator matrix is therefore exact when the dictionary is generator invariant, and in general its error is driven by the residuals $r_j=\cE_N^\perp\psi_j$, that is, by the invariance defect that also defines $\cD_N$. The residual form \eqref{eq:residual_identity} collects the statistical and regularization errors in one expression.

\begin{proposition}\label[proposition]{prop:perturbation}
Let (H1)--(H2) hold and assume the Gram smallness condition
\begin{equation}\label{eq:gram_smallness}
  c_GF_0^2\,\eta_{m,N,\delta}\le\frac{g_*}{2} .
\end{equation}
Then on an event of probability at least $1-\delta$,
\begin{equation}\label{eq:unified_bound}
  \|\widehat{\cA}_{N,m,\gamma}-\cA_N\|\le\Xi_{A,N,\delta}(m,\gamma)
  \coloneqq \frac{2\kappa_N}{g_*}\Big(c_RF_0F_r\,\eta_{m,N,\delta}+\gamma\|\boldsymbol\Lambda_N\|_2\|\mathbf A_N\|_2\Big),
\end{equation}
and, for SDMD-2 on the same event, $\|\widehat{\cB}_{N,m,\gamma}-\cB_N\|\le\Xi_{B,N,\delta}(m,\gamma)$, defined as in \eqref{eq:unified_bound} with $(c_SF_0F_s,\ \mathbf B_N)$ in place of $(c_RF_0F_r,\ \mathbf A_N)$. Here $c_G$, $c_R$, $c_S$ are absolute constants. Consequently, \cref{thm:main} applies to the empirical schemes on this event with $\xi_A=\Xi_{A,N,\delta}$ and $\xi_B=\Xi_{B,N,\delta}$.
\end{proposition}

\begin{proof}
Under (H1)--(H2), the summands of $\widehat{\mathbf G}_m-\mathbf G$ and of $\mathbf R_m$ are independent, centered, and bounded in spectral norm by $2F_0^2$ and $F_0F_r$. The matrix Bernstein inequality \cite{Tropp2012}, with Hermitian dilation for the non-Hermitian sum, gives on one event of probability at least $1-\delta$
\[
  \|\widehat{\mathbf G}_m-\mathbf G\|_2\le c_GF_0^2\,\eta_{m,N,\delta},\qquad
  \|\mathbf R_m\|_2\le c_RF_0F_r\,\eta_{m,N,\delta},
\]
and the same for the residual sum of $\widehat{\mathbf J}_m$ with $c_SF_0F_s$. By \eqref{eq:gram_smallness}, $\|\widehat{\mathbf G}_m-\mathbf G\|_2\le g_*/2$ on this event, hence $\widehat{\mathbf G}_m+\gamma\boldsymbol\Lambda_N\succeq\tfrac12g_*\mathbf I$ and $\|(\widehat{\mathbf G}_m+\gamma\boldsymbol\Lambda_N)^{-1}\|_2\le2/g_*$. Inserting these bounds into \eqref{eq:residual_identity} and converting with \eqref{eq:transport} proves \eqref{eq:unified_bound}; the bound for $\widehat{\cB}$ is identical.
\end{proof}

The bound \eqref{eq:unified_bound} is uniform over $\cV_N$ and contains the residuals of all dictionary functions through $F_r$. For a fixed observable, the sampling error of the generator matrix enters only through the invariance defect of the elements to which the matrix is applied.

\begin{proposition}\label[proposition]{prop:strong_data}
Let the hypotheses of \cref{prop:perturbation} hold with $\gamma=0$, let $v_0,\dots,v_{n-1}\in\cV_N$ be independent of the samples, and let $\dt>0$. Then, on the event of \cref{prop:perturbation} intersected with an event of probability at least $1-\delta$,
\begin{equation}\label{eq:strong_data}
  \dt\sum_{j=0}^{n-1}\big\|(\widehat{\cA}_{N,m,0}-\cA_N)v_j\big\|_{\cH}
  \le\frac{2\|\mathbf G\|_2^{1/2}F_0}{g_*\sqrt{\delta m}}\;\dt\sum_{j=0}^{n-1}\|\cE_N^\perp v_j\|_{\cH} .
\end{equation}
\end{proposition}

\begin{proof}
Let $v_j=\Phi_N\mathbf v_j$. By \eqref{eq:residual_identity} with $\gamma=0$, $(\widehat{\mathbf A}_{N,m,0}-\mathbf A_N)\mathbf v_j=\widehat{\mathbf G}_m^{-1}\mathbf R_m\mathbf v_j$, and $r(x)^\top\mathbf v_j=(\cE_N^\perp v_j)(x)$, so that $\mathbf R_m\mathbf v_j$ is the mean of the $m$ independent centered vectors $\overline{\psi(x_k)}\,(\cE_N^\perp v_j)(x_k)$ and $\E\|\mathbf R_m\mathbf v_j\|_2^2\le m^{-1}F_0^2\|\cE_N^\perp v_j\|_{\cH}^2$ by (H2). By Minkowski's inequality, $\big(\E\big[(\sum_j\|\mathbf R_m\mathbf v_j\|_2)^2\big]\big)^{1/2}\le\sum_j\big(\E\|\mathbf R_m\mathbf v_j\|_2^2\big)^{1/2}\le m^{-1/2}F_0\sum_j\|\cE_N^\perp v_j\|_{\cH}$, and by Chebyshev's inequality $\sum_j\|\mathbf R_m\mathbf v_j\|_2\le F_0(\delta m)^{-1/2}\sum_j\|\cE_N^\perp v_j\|_{\cH}$ with probability at least $1-\delta$. On the event of \cref{prop:perturbation}, $\|\widehat{\mathbf G}_m^{-1}\|_2\le2/g_*$, and $\|\Phi_N\mathbf v\|_{\cH}\le\|\mathbf G\|_2^{1/2}\|\mathbf v\|_2$ for $\mathbf v\in\C^N$ completes the proof.
\end{proof}

With $v_j\coloneqq (\cQ^{[1]}_{N,\dt})^j\Proj_Nf$, the identity \eqref{eq:telescoping} gives
\[
  \big\|(\widehat{\cQ}^{[1]})^n\Proj_Nf-(\cQ^{[1]}_{N,\dt})^n\Proj_Nf\big\|_{\cH}\le\widehat M_T\,\dt\sum_{j=0}^{n-1}\big\|(\widehat{\cA}_{N,m,0}-\cA_N)v_j\big\|_{\cH},
\]
with the bound $\widehat M_T$ for the powers of the empirical map from the proof of \cref{thm:main}, so that the data term of SDMD-1 for the observable $f$ is at most $\widehat M_T$ times the right-hand side of \eqref{eq:strong_data}. It is governed by the time integral of the invariance defect $\cE_N^\perp$ along the numerical solution, in place of the bound $F_r$ on the residuals of all dictionary functions, and it vanishes when the numerical solution stays in a generator-invariant subspace of $\cV_N$. For SDMD-R the identity $\QR_{\dt}(\widehat{\cA})-\QR_{\dt}(\cA_N)=\dt(\widehat{\cA}-\cA_N)(I+\tfrac12\dt \cA_N)+\tfrac12\dt^2\widehat{\cA}(\widehat{\cA}-\cA_N)$ gives a bound of the same form, in which $\cE_N^\perp$ is applied to the SDMD-R solution $v_j$ and to its stage values $(I+\tfrac12\dt \cA_N)v_j$, and for the exponential of \cref{rem:exponential} the argument applies to the Galerkin solution $e^{t\cA_N}\Proj_Nf$ with an integral over $[0,T]$ in place of the sum. In each case the stability of the empirical map still rests on the uniform bound \eqref{eq:unified_bound}, which controls $\widehat M_T$ when $T\,\Xi_{A,N,\delta}$ is of order one. In the experiments of \cref{sec:numerics} the empirical maps remain stable for larger perturbations, which the present analysis does not explain.

\begin{remark}\label[remark]{rem:data_variants}
If the uniform bounds of (H2) fail but $\E_\rho\|\psi\|_2^4$ and $\E_\rho[\|\psi\|_2^2\|r\|_2^2]$ are finite, as for the OU dictionary of \cref{subsec:num_ou}, Markov's inequality applied to $\E\|\widehat{\mathbf G}_m-\mathbf G\|_F^2\le m^{-1}\E_\rho\|\psi\|_2^4$ and $\E\|\mathbf R_m\|_F^2=m^{-1}\E_\rho[\|\psi\|_2^2\|r\|_2^2]$ in the Frobenius norm gives the bound \eqref{eq:unified_bound} under the condition \eqref{eq:gram_smallness} with $c_GF_0^2\,\eta_{m,N,\delta}$ and $c_RF_0F_r\,\eta_{m,N,\delta}$ replaced by $(3\E_\rho\|\psi\|_2^4/(m\delta))^{1/2}$ and $(3\E_\rho[\|\psi\|_2^2\|r\|_2^2]/(m\delta))^{1/2}$, and analogously for $\widehat{\cB}_{N,m,\gamma}$ with $s$ in place of $r$. The rate $m^{-1/2}$ is kept, and the logarithm is traded for $\delta^{-1/2}$. Samples from a stationary trajectory are dependent and require a concentration inequality for Markov chains \cite{paulin2015concentration} in place of the matrix Bernstein inequality. If the samples have the invariant law of a numerical integrator, the empirical matrices estimate the matrices of that law, and their differences from $\mathbf G$, $\mathbf H$, and $\mathbf J$ enter the perturbations $\xi_A$ and $\xi_B$. Neither case is treated here.
\end{remark}

\section{Growing dictionaries: periodic Fourier estimates}\label{sec:fourier}

\subsection{The Galerkin term in the dictionary limit}
Under the usual consistency conditions, the common term $d_N(T,f)$ of \cref{thm:main} tends to zero. The operator $e^{t\cA_N}\Proj_N$ is viewed on $\cH$ through the inclusion $\cV_N\hookrightarrow\cH$. If $\Proj_N\to I$ strongly and, for every $f$ in a core of $\cA$,
\begin{equation}\label{eq:core_consistency}
  \|\Proj_Nf-f\|_{\cH}+\|\cA_N\Proj_Nf-\cA f\|_{\cH}\to0,
\end{equation}
then the uniform resolvent bound $\|(\lambda I-\cA_N)^{-1}\|\le\lambda^{-1}$, $\lambda>0$, from Galerkin dissipativity gives resolvent convergence $(\lambda I-\cA_N)^{-1}\Proj_Nf\to(\lambda I-\cA)^{-1}f$ for all $f\in\cH$. The varying-space Trotter--Kato theorem \cite[Chap.~IX, Sec.~3]{Kato1995}, \cite{ito1998trotter} then gives $\sup_{0\le t\le T}\|e^{t\cA_N}\Proj_Nf-\cK^tf\|_{\cH}\to0$ for every $f\in\cH$, where the common exponential bound is supplied by contractivity. Rates require a concrete dictionary.

\subsection{Reversible periodic diffusions and trigonometric dictionaries}\label{subsec:periodic}
Let $\mu_0(dx)\coloneqq dx/(2\pi)^d$ on $\Tor^d$ and consider the reversible uniformly elliptic diffusion with $d\rho=\varpi\,d\mu_0$ and
\[
  \cA f\coloneqq\varpi^{-1}\nabla\cdot(a\varpi\,\nabla f)=a\Delta f+b\cdot\nabla f,\qquad b\coloneqq\nabla a+a\nabla\varpi/\varpi,
\]
where $a$ and $\varpi$ are smooth, positive, and bounded above and below; for $d=1$ this is $\cA f=af''+bf'$ with $b=a'+a\varpi'/\varpi$. Let $\cV_K$ be the trigonometric space with cutoff $K$, spanned by $e_k(x)\coloneqq e^{ik\cdot x}$, $k\in\mathbb Z^d$, $|k|_\infty\coloneqq\max_j|k_j|\le K$, of dimension $N_K\coloneqq(2K+1)^d$, and represented by the normalized real dictionary $\psi$ of the functions $1$, $\sqrt2\cos(k\cdot x)$, $\sqrt2\sin(k\cdot x)$, with $k$ running over one representative of each pair $\pm k\ne0$; for $d=1$, $\psi(x)=(1,\sqrt2\cos x,\sqrt2\sin x,\dots,\sqrt2\cos Kx,\sqrt2\sin Kx)^\top$. The Galerkin operators and $\Theta_K$ do not depend on the choice of real or exponential basis. The subscript $K$ replaces $N=N_K$ throughout, so that $\cA_K$, $\cB_K$, $\cD_K\coloneqq \cB_K-\cA_K^2$, $\Theta_K\coloneqq \|\cD_K\|$, $\cQ_{K,\dt}$, and $\widehat{\cQ}_{K,\dt,m,\gamma}$ denote the objects of \cref{sec:setting,sec:fixedN} for $\cV_K$. The quantity $\zeta_{K,m,\delta}\coloneqq \sqrt{\log(cK/\delta)/m}+\log(cK/\delta)/m$, with an absolute constant $c$, bounds the Bernstein scale $\eta_{m,N_K,\delta}$ up to a factor depending only on $d$. Constants in this section may depend on $d$. The symbols $\lesssim$ and $\gtrsim$ denote inequalities up to a constant factor independent of $K$, $\dt$, $m$, $\gamma$, $\delta$, and the observable, and $\asymp$ means that both hold.

\begin{lemma}\label[lemma]{lem:fourier_bounds}
In the setting above, $\|\psi(x)\|_2=N_K^{1/2}$ for every $x$, the Gram matrix satisfies $c_\rho\mathbf I\preceq \mathbf G_K\preceq C_\rho\mathbf I$ with $c_\rho\coloneqq\min\varpi$ and $C_\rho\coloneqq\max\varpi$, so that $\kappa_{N_K}\le(C_\rho/c_\rho)^{1/2}$, and, with constants independent of $K$,
\[
\begin{gathered}
  \|\cA_K\|\le C_AK^2,\qquad \|\cB_K\|+\Theta_K\le CK^4,\\
  \sup_x\|\cA\psi(x)\|_2\le CK^{2}N_K^{1/2},\qquad
  \sup_x\|\cA^2\psi(x)\|_2\le CK^{4}N_K^{1/2}.
\end{gathered}
\]
In particular, (H2) holds with $g_*=c_\rho$, $F_0=N_K^{1/2}$, $F_1=CK^2N_K^{1/2}$, and $F_2=CK^4N_K^{1/2}$. Moreover, there is $C_B$ depending only on $d$ and on bounds for $a$, $\nabla a$, $\varpi$, $\nabla\varpi$ such that
\begin{equation}\label{eq:energy_comparison}
  0\le\ip{v}{\cB_Kv}=\|\cA v\|_{\cH}^2\le-C_BK^2\ip{v}{\cA_Kv}\quad(v\in\cV_K),
\end{equation}
that is, $0\preceq \cB_K\preceq-C_BK^2\cA_K$.
\end{lemma}

\begin{proof}
The components of $\psi$ are orthonormal in $L^2(\mu_0)$, and $c_\rho\le\varpi\le C_\rho$ gives the bounds on $\mathbf G_K$. Each pair $\pm k\ne0$ contributes $2\cos^2(k\cdot x)+2\sin^2(k\cdot x)=2$ to $\|\psi(x)\|_2^2$, so $\|\psi(x)\|_2^2=N_K$. Trigonometric inverse inequalities give
\[
  \|\cA v\|\le CK^2\|v\|,\qquad \|\cA^2v\|\le CK^4\|v\|\qquad(v\in\cV_K),
\]
hence the bounds on $\cA_K$ and $\cB_K$, and $\Theta_K\le\|\cB_K\|+\|\cA_K\|^2$. Pointwise, $|\cA e_k|\le\|a\|_\infty|k|^2+\|b\|_\infty|k|\le CK^2$ and $|\cA^2e_k|\le CK^4$ for $|k|_\infty\le K$, and summing over the $N_K$ modes gives the bounds on the dictionary arrays. For \eqref{eq:energy_comparison}, self-adjointness gives $\ip{v}{\cB_Kv}=\|\cA v\|_{\cH}^2$. Writing $\cA v=a\Delta v+b\cdot\nabla v$ with $b$ bounded, and using the equivalence of the $\rho$- and $\mu_0$-norms and the inverse inequality $\|\Delta v\|_{L^2(\mu_0)}\le\sqrt d\,K\|\nabla v\|_{L^2(\mu_0)}$ on $\cV_K$,
\[
\begin{aligned}
  \|\cA v\|_{\cH}^2&\le C\big(\|\Delta v\|^2_{L^2(\mu_0)}+\|\nabla v\|^2_{L^2(\mu_0)}\big)\le2CdK^2\|\nabla v\|_{L^2(\mu_0)}^2,\\
  -\ip{v}{\cA_Kv}&=\int a\varpi|\nabla v|^2\,d\mu_0\ge\min(a\varpi)\|\nabla v\|^2_{L^2(\mu_0)},
\end{aligned}
\]
the identity in the second line by integration by parts.
\end{proof}

The Galerkin term of \cref{thm:main} is controlled by the following standard estimates for the periodic Fourier--Galerkin method. Let $\Proj_K^0$ be the flat Fourier truncation, the $L^2(\mu_0)$-orthogonal projector onto $\cV_K$; the equivalence of the norms of $L^2(\rho)$ and $L^2(\mu_0)$ implies that the weighted projector $\Proj_K$ is uniformly bounded on $L^2(\mu_0)$.

\begin{lemma}\label[lemma]{lem:galerkin_rate}
For every $\ell\ge0$ and $T>0$,
\begin{equation}\label{eq:sobolev_stability}
  \sup_{0\le t\le T}\|\cK^tf\|_{H^\ell(\Tor^d)}\le C_{\ell,T}\|f\|_{H^\ell(\Tor^d)}\qquad(f\in H^\ell(\Tor^d)),
\end{equation}
for every integer $0\le q\le\ell$ and $g\in H^\ell(\Tor^d)$,
\begin{equation}\label{eq:weighted_projection}
  \|(I-\Proj_K)g\|_{H^q}\le CK^{-(\ell-q)}\|g\|_{H^\ell},
\end{equation}
and for $\ell>2$ and $f\in H^\ell(\Tor^d)$,
\begin{equation}\label{eq:galerkin_rate}
  \sup_{0\le t\le T}\|e^{t\cA_K}\Proj_Kf-\cK^tf\|_{\cH}\le C_{\ell,T}\|f\|_{H^\ell}K^{-(\ell-2)} .
\end{equation}
In particular, \eqref{eq:core_consistency} holds on the core of smooth periodic functions.
\end{lemma}

\begin{proof}
Let $u(t)\coloneqq\cK^tf$ for smooth $f$, with $\cA u=a\Delta u+b\cdot\nabla u$. For an integer $q\ge0$, apply each derivative of order at most $q$ to the equation $\partial_tu=\cA u$, test with the same derivative of $u$, integrate by parts on $\Tor^d$, and sum; uniform ellipticity, bounded coefficient derivatives, and Young's inequality give
\[
  \frac12\frac{d}{dt}\|u(t)\|_{H^q}^2+c_q\|u(t)\|_{H^{q+1}}^2\le C_q\|u(t)\|_{H^q}^2,
\]
with constants depending on finitely many derivatives of $a$ and $\varpi$, and Gr\"onwall's inequality, interpolation, and density give \eqref{eq:sobolev_stability}. For \eqref{eq:weighted_projection}, write $(I-\Proj_K)g=(I-\Proj_K)(g-\Proj_K^0g)$; the inverse inequality $\|v\|_{H^q}\le CK^q\|v\|_{L^2(\mu_0)}$ on $\cV_K$ and the flat Fourier approximation bounds give
\[
  \|(I-\Proj_K)g\|_{H^q}\le\|g-\Proj_K^0g\|_{H^q}+CK^q\|g-\Proj_K^0g\|_{L^2(\mu_0)}\le CK^{-(\ell-q)}\|g\|_{H^\ell}.
\]
For \eqref{eq:galerkin_rate}, set $u_K(t)\coloneqq e^{t\cA_K}\Proj_Kf$; the Galerkin relation for $\cA_K$ gives, after adding and subtracting $\cA_K\Proj_Ku(t)$, the error equation
\[
  \frac{d}{dt}\big(u_K(t)-\Proj_Ku(t)\big)=\cA_K\big(u_K(t)-\Proj_Ku(t)\big)+\Proj_K\cA(\Proj_K-I)u(t),
\]
whose last term is bounded in $\cH$ by $C\|(I-\Proj_K)u(t)\|_{H^2}\le CK^{-(\ell-2)}\|u(t)\|_{H^\ell}$ by \eqref{eq:weighted_projection} with $q=2$, since $\cA$ is a second-order operator with bounded coefficients. Variation of constants, contractivity of $e^{t\cA_K}$, and \eqref{eq:sobolev_stability} bound $\sup_{t\le T}\|u_K(t)-\Proj_Ku(t)\|_{\cH}$ by $C_T\|f\|_{H^\ell}K^{-(\ell-2)}$, and \eqref{eq:weighted_projection} with $q=0$ bounds $\|u(t)-\Proj_Ku(t)\|_{\cH}$ by $CK^{-\ell}\|u(t)\|_{H^\ell}$. The same bound for $\Proj_K\cA(\Proj_K-I)g$ and \eqref{eq:weighted_projection} with $q=0$ for $g$ and $\cA g$, $g\in C^\infty(\Tor^d)$, verify \eqref{eq:core_consistency}.
\end{proof}

The second main result is the following theorem. As in \cref{thm:main}, perturbed operators $\widehat{\cA}$, $\widehat{\cB}$ on $\cV_K$ with $\|\widehat{\cA}-\cA_K\|\le\xi_A$ and $\|\widehat{\cB}-\cB_K\|\le\xi_B$ cover exact and empirical matrices at once.

\begin{theorem}[Periodic Fourier dictionaries]\label[theorem]{thm:fourier}
In the reversible periodic setting, let $C_A$ and $C_B$ be the constants of \cref{lem:fourier_bounds}, let $T=n\dt$, and assume the step restriction
\begin{equation}\label{eq:cfl}
  \dt\,C_AK^2\le2\qquad\text{and}\qquad\dt\,C_BK^2\le2 .
\end{equation}
Then $\cQ^{[1]}_{K,\dt}$, $\QR_{K,\dt}$, and $\cQ^{[2]}_{K,\dt}$ are self-adjoint contractions on $\cV_K$, and for $f\in H^\ell(\Tor^d)$ and perturbed operators with $\xi_A+\dt\,\xi_B\le L$,
\begin{equation}\label{eq:fourier}
\begin{aligned}
  \|(\widehat{\cQ}^{[1]})^n\Proj_Kf-\cK^Tf\|_{\cH}&\le C\|f\|_{H^\ell}\big(K^{-(\ell-2)}+\dt\big)+C\xi_A\|f\|_{\cH},\\
  \|(\widehat{\cQ}^{\mathrm R})^n\Proj_Kf-\cK^Tf\|_{\cH}&\le C\|f\|_{H^\ell}\big(K^{-(\ell-2)}+\dt^2\big)+C\xi_A\|f\|_{\cH},\\
  \|(\widehat{\cQ}^{[2]})^n\Proj_Kf-\cK^Tf\|_{\cH}&\le C\|f\|_{H^\ell}\big(K^{-(\ell-2)}+\dt^2\big)+C(\xi_A+\dt\,\xi_B)\|f\|_{\cH},
\end{aligned}
\end{equation}
where the first bound requires $\ell\ge4$ and the other two require $\ell\ge6$, with $C$ depending on $d$, $\ell$, $T$, $L$, $C_A$, $C_B$, and the coefficients, but not on $K$, $\dt$, or $f$; in particular, $\Theta_K$ does not enter the error of SDMD-2 against the Koopman semigroup on observables of this regularity. For the empirical matrices from $m$ i.i.d.\ samples under (H1) with the penalty $\boldsymbol\Lambda_{N_K}=\mathbf I$ and under the Gram smallness condition \eqref{eq:gram_smallness} with $N=N_K$, on an event of probability at least $1-\delta$ one can take
\begin{equation}\label{eq:fourier_Xi}
\begin{aligned}
  \xi_A&=\Xi_{A,K,\delta}\le C\big(K^{d+2}\zeta_{K,m,\delta}+K^2\gamma\big),\\
  \xi_B&=\Xi_{B,K,\delta}\le C\big(K^{d+4}\zeta_{K,m,\delta}+K^4\gamma\big),
\end{aligned}
\end{equation}
where $C$ is independent of $K$, $m$, $\gamma$, and $\delta$.
\end{theorem}

\begin{proof}
The three maps are real linear combinations of the self-adjoint operators $I$, $\cA_K$, $\cA_K^2$, and $\cB_K$, hence self-adjoint. By the spectral theorem, $\|\cQ^{[1]}_{K,\dt}\|=\max|1-y|$ and $\|\QR_{K,\dt}\|=\max|1-y+\tfrac12y^2|$ over $y=-\dt\lambda\in[0,\dt\|\cA_K\|]\subset[0,2]$ with $\lambda$ in the spectrum of $\cA_K$, and both maxima are at most one. For $\cQ^{[2]}_{K,\dt}=I+\dt \cA_K+\tfrac12\dt^2\cB_K$, \eqref{eq:energy_comparison} and $\cA_K\preceq0$ give
\[
  I+\dt \cA_K\preceq \cQ^{[2]}_{K,\dt}\preceq I+\dt(1-\tfrac12\dt C_BK^2)\cA_K\preceq I,
\]
and $\dt\|\cA_K\|\le2$ gives $I+\dt \cA_K\succeq-I$; hence $-I\preceq \cQ^{[2]}_{K,\dt}\preceq I$, and the three maps are contractions. This part of the proof uses \eqref{eq:cfl} only through $\dt\|\cA_K\|\le2$ and $\dt \cB_K\preceq-2\cA_K$.

For the error bounds, let first $\xi_A=\xi_B=0$, and let $u(t)\coloneqq\cK^tf$, so that $\sup_{0\le t\le T}\|u(t)\|_{H^\ell}\le C\|f\|_{H^\ell}$ by \eqref{eq:sobolev_stability}. For $\cQ$ one of the three maps set
\[
  z_j\coloneqq \cQ^j\Proj_Kf-\Proj_Ku(j\dt),\qquad
  \vartheta_j\coloneqq \cQ\Proj_Ku(j\dt)-\Proj_Ku((j+1)\dt),
\]
so that $z_0=0$ and $z_{j+1}=\cQ z_j+\vartheta_j$ with the one-step defect $\vartheta_j$. By contractivity, $\|z_n\|\le\sum_{j<n}\|\vartheta_j\|$, and $\|(I-\Proj_K)u(T)\|\le CK^{-\ell}\|u(T)\|_{H^\ell}$ by \eqref{eq:weighted_projection} accounts for the projection of the exact solution. To bound the defect, expand $u(t+\dt)$ by \eqref{eq:dynkin} with $p=1$ for SDMD-1 and $p=2$ otherwise, and write $\omega_{p+1}$ for the remainder, which satisfies $\|\omega_{p+1}\|\le C\dt^{p+1}\|\cA^{p+1}u(t)\|\le C\dt^{p+1}\|u(t)\|_{H^{2p+2}}$. Using $\cA_K\Proj_K=\Proj_K\cA\Proj_K$, $\cB_K\Proj_K=\Proj_K\cA^2\Proj_K$, and writing $u$ for $u(j\dt)$,
\begin{align*}
  \vartheta_j^{[1]}&=\dt\,\Proj_K\cA(\Proj_K-I)u-\Proj_K\omega_2,\\
  \vartheta_j^{[2]}&=\dt\,\Proj_K\cA(\Proj_K-I)u+\tfrac12\dt^2\,\Proj_K\cA^2(\Proj_K-I)u-\Proj_K\omega_3,\\
  \vartheta_j^{\mathrm R}&=\dt\,\Proj_K\cA(\Proj_K-I)u+\tfrac12\dt^2\,\Proj_K\cA\big[\Proj_K\cA(\Proj_K-I)u+(\Proj_K-I)\cA u\big]\\
  &\quad-\Proj_K\omega_3 .
\end{align*}
The projection estimate \eqref{eq:weighted_projection} gives
\[
\begin{aligned}
  \|\Proj_K\cA(\Proj_K-I)u\|&\le C\|(I-\Proj_K)u\|_{H^2}\le CK^{-(\ell-2)}\|u\|_{H^\ell},\\
  \|\Proj_K\cA^2(\Proj_K-I)u\|&\le C\|(I-\Proj_K)u\|_{H^4}\le CK^{-(\ell-4)}\|u\|_{H^\ell}.
\end{aligned}
\]
In $\vartheta_j^{\mathrm R}$, $\|\Proj_K\cA g\|\le C_AK^2\|g\|$ for $g\in\cV_K$ bounds the contribution of the first term in the bracket by $CK^2K^{-(\ell-2)}\|u\|_{H^\ell}$, and $\|\Proj_K\cA(\Proj_K-I)\cA u\|\le C\|(\Proj_K-I)\cA u\|_{H^2}\le CK^{-(\ell-4)}\|\cA u\|_{H^{\ell-2}}\le CK^{-(\ell-4)}\|u\|_{H^\ell}$ bounds that of the second. Hence
\[
  \|\vartheta_j\|\le C\|f\|_{H^\ell}\big(\dt K^{-(\ell-2)}+\dt^2K^{-(\ell-4)}+\dt^{p+1}\big),
\]
with the middle term absent for $p=1$. Summing over $n=T/\dt$ steps and using $\dt K^2\le2/C_B$ from \eqref{eq:cfl} to absorb $\dt K^{-(\ell-4)}\le CK^{-(\ell-2)}$ proves \eqref{eq:fourier} for exact matrices. For perturbed operators, split $\widehat{\cQ}^n\Proj_Kf-\cK^Tf=[\widehat{\cQ}^n-\cQ^n]\Proj_Kf+[\cQ^n\Proj_Kf-\cK^Tf]$. The powers of $\cQ$ are contractions; \cref{lem:insertion} with $\dt\|\cA_K\|\le2$, $\xi_A\le L$, and $\dt\le T$ gives $\|\widehat{\cQ}-\cQ\|\le C\dt\,\xi_A$ for SDMD-1 and SDMD-R and $\le\dt(\xi_A+\tfrac12\dt\xi_B)$ for SDMD-2, with $C$ independent of $K$; \eqref{eq:stability_transfer} with $M_T=1$ and \eqref{eq:telescoping} then bound the first term by $C_{L,T}\,\xi_A\|f\|$, respectively $C_{L,T}(\xi_A+\dt\xi_B)\|f\|$.

For \eqref{eq:fourier_Xi}, \cref{lem:fourier_bounds} gives $g_*=c_\rho$, $\kappa_{N_K}\le(C_\rho/c_\rho)^{1/2}$, $\|\mathbf A_K\|_2\lesssim K^2$, $\|\mathbf B_K\|_2\lesssim K^4$, $F_0F_r\lesssim N_K^{1/2}(K^2N_K^{1/2}+K^2N_K^{1/2})\lesssim K^2N_K\lesssim K^{d+2}$, and $F_0F_s\lesssim K^4N_K\lesssim K^{d+4}$. Inserting these into \eqref{eq:unified_bound} and its $\widehat{\cB}$ analogue, with $\log(2N_K/\delta)\le d\log(cK/\delta)$, proves \eqref{eq:fourier_Xi}.
\end{proof}

Under the step restriction \eqref{eq:cfl}, the additional data term of SDMD-2 satisfies $\dt\,\Xi_{B,K,\delta}\lesssim K^{d+2}\zeta_{K,m,\delta}+K^2\gamma$, which is the scale of $\Xi_{A,K,\delta}$. The extra matrix therefore adds to the cost but does not change the sampling rate. The sufficient condition $\Theta_K\dt\to0$ that \cref{thm:main} would impose on SDMD-2 relative to the Galerkin semigroup is not needed against the Koopman semigroup on $H^\ell$ observables, $\ell\ge6$. \Cref{ex:theta_growth} ($d=1$) shows that $\Theta_K\dt$ can grow like $K^2$ under \eqref{eq:cfl} while the error stays $O(K^{-(\ell-2)}+\dt^2)$, because the defect sits on the most strongly damped modes. For the exponential of \cref{rem:exponential}, \eqref{eq:exponential} and \eqref{eq:galerkin_rate} give the error bound $C\|f\|_{H^\ell}K^{-(\ell-2)}+T\xi_Ae^{T\xi_A}\|f\|_{\cH}$ for $\ell>2$, without a step restriction. With the largest step that \eqref{eq:cfl} allows, $\dt\asymp K^{-2}$, the time-stepping term $\dt^2\asymp K^{-4}$ of the second-order schemes is of the order of the Galerkin term for $\ell=6$ and exceeds it for smoother observables, for which a smaller step or the exponential is more accurate.

\begin{corollary}\label[corollary]{cor:accuracy_cost}
Consider the empirical matrices of \cref{thm:fourier} and $f\in H^\ell(\Tor^d)$. Fix $\epsilon>0$, $0<\delta<1/2$, and positive constants $C_m$, $c_\gamma$, $c_t$, set $\nu\coloneqq\ell-2$, $K_\epsilon\coloneqq\lceil\epsilon^{-1/\nu}\rceil$, choose $c_{\mathrm{stab}}>0$ with $c_{\mathrm{stab}}\max(C_A,C_B)\le2$, and take
\begin{equation}\label{eq:accuracy_choices}
\begin{aligned}
  m_\epsilon&\ge C_m\,\epsilon^{-2}K_\epsilon^{2d+4}\log\!\Big(\frac{cK_\epsilon}{\delta}\Big),
  &\gamma_\epsilon&\le c_\gamma\,\epsilon K_\epsilon^{-2},\\
  \dt^{[1]}_\epsilon&\le\min\{c_{\mathrm{stab}}K_\epsilon^{-2},c_t\epsilon\},
  &\dt^{\mathrm R}_\epsilon=\dt^{[2]}_\epsilon&\le\min\{c_{\mathrm{stab}}K_\epsilon^{-2},c_t\epsilon^{1/2}\},
\end{aligned}
\end{equation}
with integer step counts. Then, for sufficiently small $\epsilon$, each of the three errors is at most $C\epsilon$ with probability at least $1-\delta$, for $\ell\ge4$ in the case of SDMD-1 and $\ell\ge6$ otherwise, with $C$ independent of $\epsilon$ and $\delta$. The three schemes use the same samples, $m_\epsilon\gtrsim\epsilon^{-2-(2d+4)/\nu}\log(c\epsilon^{-1/\nu}/\delta)$, SDMD-R and SDMD-2 use the same number of steps $n_\epsilon\coloneqq T/\dt^{\mathrm R}_\epsilon$, which is of order $\max(K_\epsilon^2,\epsilon^{-1/2})$ for the largest steps allowed by \eqref{eq:accuracy_choices}, and SDMD-2 computes in addition $\widehat{\mathbf J}_m$ from $m_\epsilon N_{K_\epsilon}$ evaluations of $\cA^2\psi_j$ (\cref{tab:cost}). Since $m_\epsilon\gg n_\epsilon$ and $m_\epsilon\gg N_{K_\epsilon}$, the computation of the matrices dominates, and SDMD-R is the cheaper of the two second-order schemes by the ratio of the costs of their matrices.
\end{corollary}

\begin{proof}
The choices \eqref{eq:accuracy_choices} give $K_\epsilon^{-\nu}\lesssim\epsilon$, $\Xi_{A,K_\epsilon,\delta}\lesssim K_\epsilon^{d+2}\zeta_{K_\epsilon,m_\epsilon,\delta}+K_\epsilon^2\gamma_\epsilon\lesssim\epsilon$ by \eqref{eq:fourier_Xi}, $\dt\,\Xi_{B,K_\epsilon,\delta}\lesssim K_\epsilon^{-2}(K_\epsilon^{d+4}\zeta_{K_\epsilon,m_\epsilon,\delta}+K_\epsilon^4\gamma_\epsilon)\lesssim\epsilon$, the step restriction \eqref{eq:cfl}, and the time-stepping terms $\dt^{[1]}_\epsilon\lesssim\epsilon$ and $(\dt^{\mathrm R}_\epsilon)^2\lesssim\epsilon$; the Gram smallness condition \eqref{eq:gram_smallness} holds for small $\epsilon$ because $F_0^2\eta_{m_\epsilon,N_{K_\epsilon},\delta}\lesssim K_\epsilon^{d}\sqrt{\log(cK_\epsilon/\delta)/m_\epsilon}\lesssim\epsilon K_\epsilon^{-2}$. \Cref{thm:fourier} gives the error bounds, and the cost statements are those of \cref{tab:cost}.
\end{proof}

\begin{example}\label[example]{ex:theta_growth}
For $d=1$, let $a(x)\coloneqq1+\beta\cos x$ with $0<|\beta|<1$ and $\varpi\equiv1$, so that $\cA f=\partial_x((1+\beta\cos x)\partial_xf)$ is the generator of $dX_t=a'(X_t)\,dt+\sqrt{2a(X_t)}\,dW_t$ with invariant law $\mu_0$. Its Fourier action is
\[
  \cA e_k=-k^2e_k-\frac{\beta}{2}(k^2+k)e_{k+1}-\frac{\beta}{2}(k^2-k)e_{k-1},
\]
so only the images of $e_{\pm K}$ leave $\cV_K$, $\cE_K^\perp v=-\tfrac{\beta}{2}(K^2+K)(\hat v_Ke_{K+1}+\hat v_{-K}e_{-K-1})$ for $v=\sum_{|k|\le K}\hat v_ke_k$, and \cref{prop:residual_structure} gives
\begin{equation}\label{eq:theta_growth}
  \Theta_K=\|\cE_K^\perp\|^2=\frac{\beta^2}{4}(K^2+K)^2\asymp K^4 .
\end{equation}
For the largest steps that \eqref{eq:cfl} allows, $\dt\asymp K^{-2}$, the product $\Theta_K\dt\asymp K^2$ grows without bound, while \cref{thm:fourier} gives the error $O(K^{-(\ell-2)}+\dt^2)$ for SDMD-2 as for SDMD-R. The residual $\cD_K$ acts on $e_{\pm K}$, the modes most strongly damped by $e^{t\cA_K}$, so $\cL_{K,T}$ is small although $\Theta_K$ is large; $\Theta_K$ measures the size of the residual, not its effect over a fixed horizon. By \cref{prop:strong_data}, the sampling error of the generator matrix likewise acts on an observable only through the coefficients of $e_{\pm K}$ of its numerical solution. \Cref{subsec:num_fourier,subsec:num_spline} use this model with $\beta=1/2$.
\end{example}

\section{Numerical experiments}\label{sec:numerics}
The experiments test the two main results. The first example measures the accuracy against the Koopman semigroup and the cost of the three schemes in the periodic Fourier model of \cref{ex:theta_growth} and in a two-dimensional analogue (\cref{thm:fourier,cor:accuracy_cost}). The second isolates the time-stepping and data terms of \cref{thm:main} for a fixed non-invariant dictionary of an OU process, where the residual sits on slowly decaying modes. The third uses a spline dictionary that lies in $\dom(\cA)$ but not in $\dom(\cA^2)$. Next to the three schemes, the tables report the exponential of the generator matrix, $e^{T\cA_N}$ for exact and $e^{T\widehat{\cA}}$ for empirical matrices, which is the limit of the schemes as $\dt\to0$ (\cref{rem:exponential}). Errors against $\cK^Tf$ are measured in $L^2(\rho)$; operator errors against $e^{T\cA_N}$ are measured in the $\cV_N\to\cV_N$ norm of \eqref{eq:transport}. The horizon is $T=1$ throughout, empirical matrices are computed with $\gamma=0$, with $\widehat{\mathbf A}$, $\widehat{\mathbf B}$, $\widehat{\cA}$, $\widehat{\cB}$ written for $\widehat{\mathbf A}_{N,m,0}$, $\widehat{\mathbf B}_{N,m,0}$, $\widehat{\cA}_{N,m,0}$, $\widehat{\cB}_{N,m,0}$, the sampled computations use the seed $20260507$ ($20261029$ in the two-dimensional model), and the computations with exact matrices are deterministic. The code that reproduces every table and figure is available on \href{https://github.com/TalkingDoll/sdmd_analysis.git}{GitHub}.

\subsection{Periodic Fourier models: accuracy and cost}\label{subsec:num_fourier}
For the model of \cref{ex:theta_growth} with $\beta=1/2$, the Galerkin matrices are known in closed form in the orthonormal basis $\{e_k\}$, tridiagonal for $\cA_K$ and pentadiagonal for $\cB_K$, so the errors with exact matrices require no quadrature. The observable has Fourier coefficients $\hat f_k\propto(1+k^2)^{-7/2}$, normalized to $\|f\|=1$, and lies in $H^\ell(\Tor)$ for every $\ell<6.5$; the reference $\cK^Tf$ is the Galerkin solution with cutoff $256$, computed with the matrix exponential. For each cutoff the step is $\dt=T/n$ with $n$ the smallest integer such that $\dt\|\cA_K\|\le1.5$. The conditions $\dt\|\cA_K\|\le2$ and $\dt \cB_K\preceq-2\cA_K$, through which the proof of \cref{thm:fourier} uses the step restriction \eqref{eq:cfl}, hold in the computation for every $K$ of \cref{tab:fourier_population}, and the computed norm of $\cQ^{[2]}_{K,\dt}$ equals one.

\Cref{tab:fourier_population} reports the errors with exact matrices. SDMD-R and SDMD-2 agree to three digits for $K\ge8$ and decay like $K^{-4}\asymp\dt^2$ (fitted slopes $-4.15$ over $K=16$--$64$), SDMD-1 decays like $K^{-2}\asymp\dt$ (slope $-2.07$), and at $K=4$ SDMD-R is more accurate than SDMD-2 by five percent. The Galerkin error, which is the error of the exponential $e^{T\cA_K}$, is smaller than the error of SDMD-R by a factor that grows from $1.4$ at $K=4$ to $10^{30}$ at $K=64$, because $\cK^Tf$ is analytic; the bound $K^{-(\ell-2)}$ of \eqref{eq:galerkin_rate} is uniform in $t\in[0,T]$ and is pessimistic at $t=T$. For $K\ge8$ the errors of the three schemes are therefore time-stepping errors, in agreement with \eqref{eq:fourier}, although $\Theta_K\dt$ grows from $2.1$ at $K=4$ to $291$ at $K=64$; the closed form \eqref{eq:theta_growth} is reproduced to machine precision for $1\le K\le64$.

\begin{table}[tbp]
\footnotesize
\centering
\setlength{\tabcolsep}{5pt}
\begin{tabular}{rrrrrrr}
\toprule
$K$ & $n$ & $\Theta_K\dt$ & Exponential & SDMD-1 & SDMD-R & SDMD-2 \\
\midrule
4 & 12 & $2.1$ & $7.4\times10^{-5}$ & $1.78\times10^{-3}$ & $1.03\times10^{-4}$ & $1.08\times10^{-4}$ \\
8 & 49 & $6.6$ & $2.1\times10^{-7}$ & $4.26\times10^{-4}$ & $3.61\times10^{-6}$ & $3.61\times10^{-6}$ \\
16 & 210 & $22$ & $3.0\times10^{-12}$ & $9.89\times10^{-5}$ & $1.91\times10^{-7}$ & $1.91\times10^{-7}$ \\
32 & 889 & $78$ & $1.1\times10^{-21}$ & $2.33\times10^{-5}$ & $1.06\times10^{-8}$ & $1.06\times10^{-8}$ \\
64 & 3716 & $291$ & $2.8\times10^{-40}$ & $5.58\times10^{-6}$ & $6.05\times10^{-10}$ & $6.05\times10^{-10}$ \\
\bottomrule
\end{tabular}
\captionbelow{Fourier model with exact Galerkin matrices, $\beta=1/2$, $T=1$, and $\dt=T/n$ with $\dt\|\cA_K\|\le1.5$: errors $\|\cP\,\Proj_Kf-\cK^Tf\|_{\cH}$ for the exponential $\cP=e^{T\cA_K}$, which is the Galerkin error, and for $\cP=(\cQ_{K,\dt})^n$ with the three schemes. The Galerkin error is evaluated in $60$-digit arithmetic with reference cutoff $96$, whose truncation error is below $10^{-58}$. Fitted slopes in $K$ over $K=16$--$64$: $-2.07$ for SDMD-1 and $-4.15$ for SDMD-R and SDMD-2.}
\label{tab:fourier_population}
\end{table}

\Cref{tab:fourier_sampled} repeats the computation with empirical matrices from $m$ i.i.d.\ samples of the uniform invariant law, computed by the generic formulas $\cA\psi=a\psi''+a'\psi'$ and $\cA^2\psi=\cA(\cA\psi)$ from the derivatives of the dictionary and of $a$. Although the matrix perturbations are of order one and $10^2$ at $K=8$, the mean errors of the three schemes differ from the errors with exact matrices of \cref{tab:fourier_population} by at most $0.3$ percent for every $(K,m)$; those of SDMD-1 agree to three digits and are not listed. This agreement is explained by the residual identity \eqref{eq:residual_identity} and \cref{prop:strong_data}. The residuals of the dictionary functions with wave numbers below $K$ vanish because $\cA e_k\in\cV_K$ for $|k|<K$, so $\widehat{\mathbf A}-\mathbf A_K=\widehat{\mathbf G}_m^{-1}\mathbf R_m$ is supported on the two columns of the wave number $K$, and the sampling error acts only on the coefficients of $e_{\pm K}$ of the numerical solution, which start from $\hat f_{\pm K}\propto(1+K^2)^{-7/2}$ and decay at the fastest rate; likewise the residuals of $\cA^2$ vanish for wave numbers below $K-1$. The last two columns of \cref{tab:fourier_sampled} are the two sides of \eqref{eq:strong_data} for the SDMD-1 solution $v_j=(\cQ^{[1]}_{K,\dt})^j\Proj_Kf$ with exact matrices, where $\mathbf G=\mathbf I$, $g_*=1$, $F_0=N_K^{1/2}$, and $\delta=0.1$. The left-hand side is about one sixth of the right-hand side in every row and exceeds it in no repetition. For $m=10^4$ it decreases from $3.7\times10^{-5}$ at $K=4$ to $1.9\times10^{-14}$ at $K=64$, while $\|\widehat{\mathbf A}-\mathbf A_K\|_2$ grows from $0.17$ to $1.3\times10^{2}$. A data term of the same form enters the exponential $e^{T\widehat{\cA}}$ of the empirical generator matrix, whose error exceeds the Galerkin error of \cref{tab:fourier_population} by $2\times10^{-6}$ at $K=4$ and by less than $10^{-11}$ for $K\ge16$. The perturbations grow like $K^{2.4}$ and $K^{4.3}$ at fixed $m$ (fitted over $K=4$--$64$), below the bounds $K^3$ and $K^5$ of \eqref{eq:fourier_Xi}, and decay like $m^{-1/2}$ at fixed $K$.

\begin{table}[tbp]
\ifjsc\scriptsize\else\footnotesize\fi
\centering
\setlength{\tabcolsep}{\ifjsc 2pt\else 4pt\fi}
\begin{tabular}{rrrrrrrrr}
\toprule
 & & & & & \multicolumn{2}{c}{Perturbation of} & \multicolumn{2}{c}{Sides of \eqref{eq:strong_data}} \\
\cmidrule(lr){6-7}\cmidrule(lr){8-9}
$K$ & $m$ & Exponential & SDMD-R & SDMD-2 & $\mathbf A_K$ & $\mathbf B_K$ & left & right \\
\midrule
4 & $10^{4}$ & $7.6\times10^{-5}$ & $1.03\times10^{-4}$ & $1.08\times10^{-4}$ & $0.17$ & $7.2$ & $3.7\times10^{-5}$ & $2.3\times10^{-4}$ \\
8 & $10^{4}$ & $2.4\times10^{-7}$ & $3.60\times10^{-6}$ & $3.61\times10^{-6}$ & $0.88$ & $1.3\times10^{2}$ & $3.6\times10^{-7}$ & $2.1\times10^{-6}$ \\
16 & $10^{4}$ & $7.8\times10^{-12}$ & $1.91\times10^{-7}$ & $1.91\times10^{-7}$ & $4.4$ & $2.4\times10^{3}$ & $2.7\times10^{-10}$ & $1.7\times10^{-9}$ \\
32 & $10^{4}$ & $<10^{-12}$ & $1.06\times10^{-8}$ & $1.06\times10^{-8}$ & $23$ & $5.1\times10^{4}$ & $2.1\times10^{-12}$ & $1.3\times10^{-11}$ \\
64 & $10^{4}$ & $<10^{-12}$ & $6.05\times10^{-10}$ & $6.05\times10^{-10}$ & $1.3\times10^{2}$ & $1.1\times10^{6}$ & $1.9\times10^{-14}$ & $1.2\times10^{-13}$ \\
\midrule
8 & $10^{3}$ & $3.8\times10^{-7}$ & $3.60\times10^{-6}$ & $3.60\times10^{-6}$ & $2.8$ & $4.0\times10^{2}$ & $1.2\times10^{-6}$ & $6.7\times10^{-6}$ \\
8 & $10^{5}$ & $2.2\times10^{-7}$ & $3.61\times10^{-6}$ & $3.61\times10^{-6}$ & $0.25$ & $37$ & $1.0\times10^{-7}$ & $6.7\times10^{-7}$ \\
32 & $10^{5}$ & $<10^{-12}$ & $1.06\times10^{-8}$ & $1.06\times10^{-8}$ & $7.3$ & $1.6\times10^{4}$ & $6.7\times10^{-13}$ & $4.2\times10^{-12}$ \\
\bottomrule
\end{tabular}
\captionbelow{Fourier model with empirical matrices from $m$ i.i.d.\ uniform samples and the steps of \cref{tab:fourier_population}, means over ten repetitions: errors $\|\cP\,\Proj_Kf-\cK^Tf\|_{\cH}$ for the exponential $\cP=e^{T\widehat{\cA}}$ and for SDMD-R and SDMD-2, perturbations $\|\widehat{\mathbf A}-\mathbf A_K\|_2$ and $\|\widehat{\mathbf B}-\mathbf B_K\|_2$ of the Galerkin matrices, and the two sides of \eqref{eq:strong_data} for the SDMD-1 solution with exact matrices and $\delta=0.1$. Errors below $10^{-12}$ are not resolved in double precision.}
\label{tab:fourier_sampled}
\end{table}

The same behavior is observed in two dimensions. On $\Tor^2$ let $a(x)\coloneqq1+\tfrac14\cos((1,1)\cdot x)+\tfrac14\cos((1,-1)\cdot x)$ and $\varpi\equiv1$, so that $\cA f=\nabla\cdot(a\nabla f)$ generates $dX_t=\nabla a(X_t)\,dt+\sqrt{2a(X_t)}\,dW_t$ with uniform invariant law; $\cA e_k$ is a combination of $e_k$ and the four modes $e_{k+(\pm1,\pm1)}$, so the Galerkin matrices are again explicit. The dictionary consists of the $N_K=(2K+1)^2$ modes with $|k|_\infty\le K$, the observable has $\hat f_k\propto(1+|k|^2)^{-4}$, normalized to $\|f\|=1$ and in $H^\ell(\Tor^2)$ for $\ell<7$, and the step is chosen as before; the reference solution is the Galerkin solution with cutoff $64$, which differs from those with cutoffs $48$ and $80$ by $5\times10^{-13}$. \Cref{tab:torus2d} reports the errors with exact matrices and with $m=10^4$ samples. The three one-step maps are contractions, SDMD-1 converges at first order and the two second-order schemes at second order in $\dt\asymp K^{-2}$ (fitted slopes $-2.1$ and $-4.2$ in $K$ over $K=8$--$32$), and the errors of the two second-order schemes agree to within ten percent at $K=2$ and to within one percent for $K\ge8$. With empirical matrices the mean errors of the three schemes differ from those with exact matrices by at most five percent, also at $K=16$, where $m/N_K\approx9$ and $\|\widehat{\mathbf A}-\mathbf A_K\|_2\approx54$.

\begin{table}[tbp]
\footnotesize
\centering
\setlength{\tabcolsep}{5pt}
\begin{tabular}{rrrrrr}
\toprule
$K$ & $n$ & Exponential & SDMD-1 & SDMD-R & SDMD-2 \\
\midrule
\multicolumn{6}{l}{Exact matrices} \\[2pt]
2 & 6 & $6.5\times10^{-4}$ & $4.20\times10^{-3}$ & $7.47\times10^{-4}$ & $6.82\times10^{-4}$ \\
4 & 22 & $2.4\times10^{-5}$ & $1.06\times10^{-3}$ & $3.05\times10^{-5}$ & $2.97\times10^{-5}$ \\
8 & 91 & $5.9\times10^{-8}$ & $2.53\times10^{-4}$ & $1.03\times10^{-6}$ & $1.03\times10^{-6}$ \\
16 & 391 & $<10^{-12}$ & $5.87\times10^{-5}$ & $5.52\times10^{-8}$ & $5.52\times10^{-8}$ \\
32 & 1692 & $<10^{-12}$ & $1.36\times10^{-5}$ & $2.94\times10^{-9}$ & $2.94\times10^{-9}$ \\
\midrule
\multicolumn{6}{l}{Empirical matrices, $m=10^4$} \\[2pt]
2 & 6 & $6.6\times10^{-4}$ & $4.23\times10^{-3}$ & $7.43\times10^{-4}$ & $6.85\times10^{-4}$ \\
4 & 22 & $2.6\times10^{-5}$ & $1.06\times10^{-3}$ & $3.18\times10^{-5}$ & $3.00\times10^{-5}$ \\
8 & 91 & $9.6\times10^{-8}$ & $2.53\times10^{-4}$ & $1.04\times10^{-6}$ & $1.03\times10^{-6}$ \\
16 & 391 & $3.2\times10^{-12}$ & $5.87\times10^{-5}$ & $5.52\times10^{-8}$ & $5.52\times10^{-8}$ \\
\bottomrule
\end{tabular}
\captionbelow{Two-dimensional Fourier model of \cref{subsec:num_fourier}, $T=1$, and $\dt=T/n$ with $\dt\|\cA_K\|\le1.5$: errors $\|\cP\,\Proj_Kf-\cK^Tf\|_{\cH}$ for the exponential and the three schemes, with exact matrices and with empirical matrices from $m=10^4$ uniform samples (means over ten repetitions, five for $K=16$). Errors below $10^{-12}$ are not resolved by the reference solution. Fitted slopes in $K$ over $K=8$--$32$ for exact matrices: $-2.11$ for SDMD-1 and $-4.23$ for SDMD-R and SDMD-2.}
\label{tab:torus2d}
\end{table}

\begin{figure}[tbp]
\centering
\includegraphics[width=\figwidth]{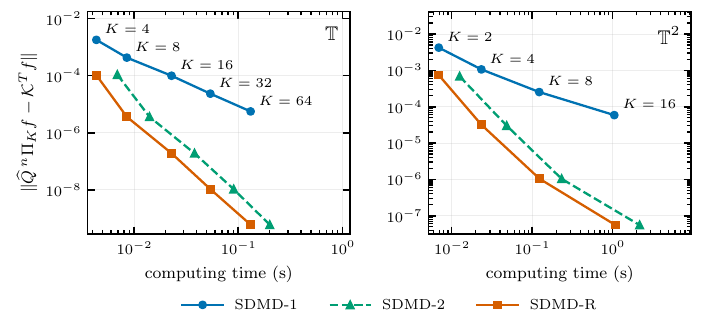}
\caption{Error against $\cK^Tf$ versus the computing time for the matrices and the time stepping, with $m=10^4$ samples, for the Fourier model on $\Tor$ with $K=4,8,16,32,64$ (left, \cref{tab:fourier_sampled,tab:times}) and on $\Tor^2$ with $K=2,4,8,16$ (right, \cref{tab:torus2d,tab:times}). SDMD-R lies below SDMD-1 at equal computing time and to the left of SDMD-2 at equal error\ifjsc\else.\fi}
\label{fig:work_precision}
\end{figure}

\Cref{tab:times} lists the computing times for $m=10^4$ samples. On $\Tor$ the matrices of SDMD-2 take $1.6$ to $1.7$ times as long to compute as those of SDMD-1 and SDMD-R, which are the same, and on $\Tor^2$ they take $1.8$ to $2.2$ times as long, so at equal cutoff and sample size SDMD-R reaches the accuracy of SDMD-2 at the cost of SDMD-1 (\cref{fig:work_precision}). On $\Tor$ the time stepping takes at most $17$~ms, the same for the three schemes up to the product $\widehat{\mathbf A}^2$ of SDMD-R, and the exponential at most $8$~ms; on $\Tor^2$ at $K=16$ they take $0.17$~s and $0.66$~s.

\begin{table}[tbp]
\footnotesize
\centering
\begin{tabular}{lrrrrrrr}
\toprule
 & $K$ & $N_K$ & $t_{1,\mathrm R}$ & $t_2$ & $t_2/t_{1,\mathrm R}$ & Time stepping & Exponential \\
\midrule
$\Tor$ & 4 & 9 & $0.0043$ & $0.0069$ & $1.6$ & $1.5\times10^{-5}$ & $1.9\times10^{-5}$ \\
 & 8 & 17 & $0.0085$ & $0.014$ & $1.7$ & $4.4\times10^{-5}$ & $4.4\times10^{-5}$ \\
 & 16 & 33 & $0.023$ & $0.038$ & $1.7$ & $2.0\times10^{-4}$ & $2.2\times10^{-4}$ \\
 & 32 & 65 & $0.052$ & $0.089$ & $1.7$ & $0.0017$ & $0.0011$ \\
 & 64 & 129 & $0.11$ & $0.18$ & $1.6$ & $0.017$ & $0.0071$ \\
\midrule
$\Tor^2$ & 2 & 25 & $0.0069$ & $0.013$ & $1.8$ & $1.1\times10^{-5}$ & $3.0\times10^{-5}$ \\
 & 4 & 81 & $0.023$ & $0.048$ & $2.1$ & $6.1\times10^{-5}$ & $3.9\times10^{-4}$ \\
 & 8 & 289 & $0.12$ & $0.23$ & $1.9$ & $0.0019$ & $0.015$ \\
 & 16 & 1089 & $0.92$ & $2.0$ & $2.2$ & $0.17$ & $0.66$ \\
\bottomrule
\end{tabular}
\captionbelow{Computing times in seconds for $m=10^4$ samples in the Fourier models on $\Tor$ and $\Tor^2$, best of three runs on a two-core Linux machine with OpenBLAS: matrices of SDMD-1 and SDMD-R ($t_{1,\mathrm R}$, for the arrays $\boldsymbol\Psi$, $\cA\boldsymbol\Psi$ and the matrices $\widehat{\mathbf G}_m$, $\widehat{\mathbf H}_m$, $\widehat{\mathbf A}$), matrices of SDMD-2 ($t_2$, in addition $\cA^2\boldsymbol\Psi$, $\widehat{\mathbf J}_m$, $\widehat{\mathbf B}$), time stepping of SDMD-R (the product $\widehat{\mathbf A}^2$ and $n$ matrix--vector products), and the exponential $e^{T\widehat{\mathbf A}}$.}
\label{tab:times}
\end{table}

In both Fourier models $\cK^Tf$ is analytic, and the Galerkin error decays faster in $K$ than the time-stepping errors at the steps $\dt\asymp K^{-2}$ used here. For $K\ge8$ it is more than one order of magnitude below the errors of the three schemes. The exponential of the generator matrix is the most accurate of the propagators in every row of \cref{tab:fourier_population,tab:fourier_sampled,tab:torus2d}, with empirical matrices by one order of magnitude at $K=8$ and by more than four at $K=16$, and its computing time is below that of the matrices in every row of \cref{tab:times}. The comparison of the three schemes at a step close to the stability limit measures their orders in $\dt$ and the effect of the third matrix. The examples of \cref{subsec:num_ou,subsec:num_spline}, in which the Galerkin error dominates, show when SDMD-R reaches the accuracy of the exponential.

\subsection{OU process with a fixed non-invariant dictionary}\label{subsec:num_ou}
The second example is the one-dimensional OU process $dX_t=-X_t\,dt+\sqrt2\,dW_t$ with $\rho=\mathcal N(0,1)$ and the dictionary $\{1,\sin x,\cos x,\sin2x,\cos2x\}$, which is not generator invariant, since $\cA\sin x=-x\cos x-\sin x\notin\cV_N$. Exact matrices are computed with $160$ Gauss--Hermite nodes. Here $\Theta_N\approx3.66$, the eigenvalues of $\cA_N$ are $0$, $-1.20$, $-2.09$, $-4.93$, $-5.44$, the leading operator $\cL_{N,T}$ of \eqref{eq:leading_term} has norm $0.108$, and the lower bound \eqref{eq:lower_bound} equals $0.082$. The invariance defect sits on slowly decaying modes, so this is the case in which the term $\Theta_N\dt$ of \cref{thm:main} is visible.

\Cref{tab:ou_population} isolates the time-stepping error, the error against the common limit $e^{T\cA_N}$ to which the three schemes converge as $\dt\to0$. SDMD-R converges at second order (fitted slope $2.01$ over the three finest steps), and SDMD-1 and SDMD-2 at first order ($1.00$ and $1.02$). SDMD-2 converges with the predicted leading term, since the distance of $((\cQ^{[2]}_{N,\dt})^n-e^{T\cA_N})/\dt$ from $\cL_{N,T}$, divided by $\|\cL_{N,T}\|$, is about halved with each halving of $\dt$, as the $O(\dt^2)$ remainder in \eqref{eq:leading_term} requires. At $\dt=0.1$ the time-stepping error of SDMD-R is seven times smaller than that of SDMD-2, at $\dt=0.003125$ nearly two hundred times. Against the Koopman semigroup itself, for $f=\cos x$ with $\cK^Tf=e^{-(1-e^{-2T})/2}\cos(e^{-T}x)$, all three schemes approach the common Galerkin error $d_N(T,f)=9.55\times10^{-3}$ of the five-dimensional dictionary as $\dt\to0$. Here SDMD-2 is the more accurate of the two second-order schemes at every step, as \cref{rem:reading} allows, but the difference is $4.2$ percent of $d_N(T,f)$ at $\dt=0.1$ and $0.1$ percent at $\dt=0.003125$; the scheme-dependent part of the error is the time-stepping part, in which SDMD-R has the higher order.

\begin{table}[tbp]
\footnotesize
\centering
\setlength{\tabcolsep}{3pt}
\begin{tabular}{rrrrrrrr}
\toprule
& \multicolumn{4}{c}{Error against $e^{T\cA_N}$} & \multicolumn{3}{c}{Error against $\cK^T\cos x$} \\
\cmidrule(lr){2-5}\cmidrule(lr){6-8}
$\dt$ & SDMD-1 & SDMD-R & SDMD-2 & rel. & SDMD-1 & SDMD-R & SDMD-2 \\
\midrule
$0.1$ & $2.78\times10^{-2}$ & $2.40\times10^{-3}$ & $1.62\times10^{-2}$ & $0.574$ & $17.5$ & $9.37$ & $8.97$ \\
$0.05$ & $1.37\times10^{-2}$ & $5.10\times10^{-4}$ & $6.47\times10^{-3}$ & $0.224$ & $12.6$ & $9.50$ & $9.30$ \\
$0.025$ & $6.80\times10^{-3}$ & $1.22\times10^{-4}$ & $2.94\times10^{-3}$ & $0.102$ & $10.8$ & $9.54$ & $9.44$ \\
$0.0125$ & $3.39\times10^{-3}$ & $3.00\times10^{-5}$ & $1.41\times10^{-3}$ & $0.049$ & $10.1$ & $9.55$ & $9.50$ \\
$0.00625$ & $1.69\times10^{-3}$ & $7.42\times10^{-6}$ & $6.88\times10^{-4}$ & $0.024$ & $9.78$ & $9.55$ & $9.53$ \\
$0.003125$ & $8.45\times10^{-4}$ & $1.85\times10^{-6}$ & $3.40\times10^{-4}$ & $0.012$ & $9.66$ & $9.55$ & $9.54$ \\
\bottomrule
\end{tabular}
\captionbelow{OU process with the trigonometric dictionary and exact Galerkin matrices, $T=1$, $n=T/\dt$: the operator norms $\|(\cQ_{N,\dt})^n-e^{T\cA_N}\|$ (fitted slopes over the three finest steps $1.00$, $2.01$, $1.02$), the normalized distance of $((\cQ^{[2]}_{N,\dt})^n-e^{T\cA_N})/\dt$ from $\cL_{N,T}$ (rel.), and $\|(\cQ_{N,\dt})^n\Proj_N\cos x-\cK^T\cos x\|_{L^2(\rho)}$ in units of $10^{-3}$, whose common limit is the Galerkin error $9.55$.}
\label{tab:ou_population}
\end{table}

\begin{table}[tbp]
\footnotesize
\centering
\begin{tabular}{rllll}
\toprule
$m$ & Exponential & SDMD-1 & SDMD-R & SDMD-2 \\
\midrule
\multicolumn{5}{l}{Error against $e^{T\cA_N}$} \\[2pt]
400 & $4.98$ [$2.4$, $8.8$] & $5.17$ [$2.5$, $8.8$] & $4.98$ [$2.4$, $8.8$] & $4.54$ [$2.1$, $8.2$] \\
800 & $3.31$ [$1.3$, $5.9$] & $3.64$ [$1.8$, $6.2$] & $3.31$ [$1.3$, $5.9$] & $3.12$ [$1.4$, $5.5$] \\
1600 & $2.25$ [$0.98$, $3.9$] & $2.67$ [$1.7$, $4.1$] & $2.25$ [$0.97$, $3.9$] & $2.17$ [$0.92$, $3.8$] \\
3200 & $1.58$ [$0.80$, $2.6$] & $2.09$ [$1.4$, $2.9$] & $1.58$ [$0.79$, $2.6$] & $1.61$ [$0.80$, $2.6$] \\
6400 & $1.08$ [$0.50$, $1.8$] & $1.85$ [$1.3$, $2.4$] & $1.08$ [$0.49$, $1.8$] & $1.16$ [$0.69$, $1.8$] \\
12800 & $0.764$ [$0.38$, $1.2$] & $1.58$ [$1.3$, $1.9$] & $0.766$ [$0.37$, $1.2$] & $0.952$ [$0.66$, $1.3$] \\
25600 & $0.527$ [$0.26$, $0.81$] & $1.49$ [$1.3$, $1.7$] & $0.529$ [$0.27$, $0.82$] & $0.796$ [$0.55$, $1.0$] \\
\midrule
\multicolumn{5}{l}{Error against $\cK^T\cos x$} \\[2pt]
400 & $1.22$ & $1.46$ & $1.22$ & $1.16$ \\
800 & $1.03$ & $1.31$ & $1.03$ & $0.997$ \\
1600 & $0.997$ & $1.28$ & $0.992$ & $0.966$ \\
3200 & $0.977$ & $1.27$ & $0.972$ & $0.949$ \\
6400 & $0.976$ & $1.29$ & $0.970$ & $0.947$ \\
12800 & $0.961$ & $1.26$ & $0.956$ & $0.935$ \\
25600 & $0.960$ & $1.27$ & $0.954$ & $0.934$ \\
\bottomrule
\end{tabular}
\captionbelow{OU process with the trigonometric dictionary and empirical matrices from $m$ i.i.d.\ samples of $\mathcal N(0,1)$, $\dt=0.05$, $100$ repetitions, all entries in units of $10^{-2}$. Top: operator errors $\|\cP-e^{T\cA_N}\|$ for the exponential $\cP=e^{T\widehat{\cA}}$ and for $\cP=(\widehat{\cQ}_{N,\dt,m,0})^n$, means with the 10th and 90th percentiles in brackets; fitted slopes in $m$ over the five smallest sample sizes $-0.55$, $-0.38$ (SDMD-1, approaching its time-stepping error), $-0.55$, $-0.49$. Bottom: mean errors $\|\cP\,\Proj_N\cos x-\cK^T\cos x\|_{L^2(\rho)}$; the Galerkin error is $0.955$.}
\label{tab:ou_sampled}
\end{table}

\Cref{tab:ou_sampled} adds the data term. With $\dt=0.05$ and empirical matrices from $m$ i.i.d.\ samples of $\rho$, the error against $e^{T\cA_N}$ of each scheme decreases like $m^{-1/2}$ until it reaches the time-stepping error of the scheme, which is $1.37\times10^{-2}$ for SDMD-1, $5.10\times10^{-4}$ for SDMD-R, and $6.47\times10^{-3}$ for SDMD-2 by \cref{tab:ou_population}. The mean error of SDMD-1 is $15$ percent above its time-stepping error at $m=12800$ and nine percent at $m=25600$, where that of SDMD-2 is $23$ percent above its own. The time-stepping error of SDMD-R is below the data term at every sample size, and the error of SDMD-R agrees with that of the exponential $e^{T\widehat{\cA}}$ to within half a percent. For $m\le1600$ the mean error of SDMD-2 is the smaller of the two, by up to nine percent. From $m=3200$ on SDMD-R is more accurate with the same samples and without the matrix $\widehat{\mathbf J}_m$, at $m=25600$ by one third and in $94$ of the $100$ repetitions. The perturbation of that matrix, $\|\widehat{\cB}-\cB_N\|$, is between $10$ and $11$ times $\|\widehat{\cA}-\cA_N\|$, which decreases from $0.66$ at $m=400$ to $0.073$ at $m=25600$, because the array $\cA^2\boldsymbol\Psi$ contains the features $x^2\cos x$ and $x\sin2x$. Against $\cK^T\cos x$, the errors of the exponential, SDMD-R, and SDMD-2 are within five percent of the Galerkin error from $m=1600$ on, and the error of SDMD-1 stays one third above it. SDMD-R and the exponential agree to within one percent, and SDMD-2 is more accurate than SDMD-R in every repetition, in the mean by two to seven percent of the Galerkin error. The dictionary violates (H2), and the rate $m^{-1/2}$ is the one of \cref{rem:data_variants}.

\subsection{A spline dictionary with two derivatives}\label{subsec:num_spline}
The third example returns to the model of \cref{ex:theta_growth} with $\beta=1/2$ and the observable of \cref{subsec:num_fourier}, with the dictionary of the $N$ periodic cubic B-splines $\psi_j(x)\coloneqq\sum_{l\in\mathbb Z}\chi((x-jh)/h-lN)$, $0\le j<N$, on the uniform mesh of width $h\coloneqq2\pi/N$, where $\chi$ is the cubic cardinal B-spline with support $[-2,2]$. These functions are twice continuously differentiable and piecewise cubic. They lie in $\dom(\cA)=H^2(\Tor)$, but the derivative of $\cA\psi_j$ jumps at the mesh points, so $\cV_N\not\subset\dom(\cA^2)$, and the trigonometric estimates of \cref{sec:fourier} do not apply. SDMD-1 and SDMD-R are defined as before and are contractions for $\dt\|\cA_N\|\le2$, because $\cA_N$ is self-adjoint. For SDMD-2 we compare two substitutes for the matrix $\mathbf J$. The first evaluates $\cA^2\psi_j$ at the samples by the differentiation formulas; the samples almost surely avoid the mesh points, where the fourth derivative of the spline is concentrated. The second is the energy form $\mathbf J^{\mathrm{en}}_{ij}\coloneqq\ip{\cA\psi_i}{\cA\psi_j}$ of \cref{app:energy}, which defines $\cB_N$ on this dictionary. Integration by parts on each mesh cell gives
\[
  \ip{\cA\psi_i}{\cA\psi_j}=\int_{\Tor}\psi_i\,(\cA^2\psi_j)_{\mathrm{pw}}\,d\rho+\frac1{2\pi}\sum_{k=0}^{N-1}\psi_i(kh)\,a(kh)^2\,[\psi_j'''](kh),
\]
where $(\cA^2\psi_j)_{\mathrm{pw}}$ is the pointwise evaluation away from the mesh points and $[\psi_j''']$ is the jump of the third derivative. The second term, which no sample sees, is as large as the left-hand side. Exact matrices are computed by Gauss--Legendre quadrature with ten points per mesh cell, and the step is $\dt=T/n$ with $n$ the smallest integer such that $\dt\|\cA_N\|\le1.5$.

\Cref{tab:spline} reports the errors against $\cK^Tf$. With exact matrices the Galerkin error, the error of $e^{T\cA_N}$, decays like $N^{-4}$ (fitted slope $-4.10$ over $N=32$--$128$). Since $\dt\asymp\|\cA_N\|^{-1}\asymp N^{-2}$, the error of SDMD-R decays at the same rate, and it exceeds the Galerkin error by $7$ to $11$ percent. The Galerkin error is not negligible here, and SDMD-R has the accuracy of the exponential without a function of the matrix. SDMD-2 with the energy form agrees with SDMD-R to within one percent, although $\Theta_N\dt$ grows from $6.5$ at $N=16$ to $299$ at $N=128$. SDMD-2 with the pointwise evaluation is first-order accurate only, with an error $1.4$ times that of SDMD-1 (slope $-2.05$ for both). When the dictionary is not in $\dom(\cA^2)$, the sampled values of $\cA^2\psi_j$ do not determine the energy form $\ip{\cA\psi_i}{\cA\psi_j}$, and this error does not decrease with the sample size.

With empirical matrices from $m$ uniform samples, the mean error of SDMD-R stays within eight percent of that of the exponential $e^{T\widehat{\cA}}$ and approaches the error with exact matrices as $m$ grows. The energy form needs more samples. At $N=32$ the mean error of SDMD-2 with the energy form is $32$ times that of SDMD-R for $m=10^4$, $11$ times for $m=10^5$, and $3.0$ times for $m=10^6$; at $N=16$ the factors are $4.7$, $1.5$, and $1.04$. SDMD-R is the more accurate of the two in every repetition except one of the $50$ repetitions at $N=16$ and $m=10^6$. The matrix $\widehat{\mathbf J}^{\mathrm{en}}_m$ is not of residual form (\cref{app:energy}), so its sampling error acts on every element of $\cV_N$ and not only through the invariance defect, and it enters \eqref{eq:main} through the term $\dt\,\xi_B$, which SDMD-R does not have. On this dictionary SDMD-R is second-order accurate with the pointwise evaluations and the two matrices of SDMD-1, whereas SDMD-2 requires the energy form, a third matrix, and more samples.

\begin{table}[tbp]
\footnotesize
\centering
\setlength{\tabcolsep}{3.5pt}
\begin{tabular}{rrrrrrrr}
\toprule
 & & & & & & \multicolumn{2}{c}{SDMD-2} \\
\cmidrule(lr){7-8}
$N$ & $m$ & $n$ & Exponential & SDMD-1 & SDMD-R & pointwise & energy \\
\midrule
16 & exact & 54 & $7.87\times10^{-6}$ & $3.87\times10^{-4}$ & $8.45\times10^{-6}$ & $5.43\times10^{-4}$ & $8.37\times10^{-6}$ \\
32 & exact & 233 & $3.22\times10^{-7}$ & $8.91\times10^{-5}$ & $3.58\times10^{-7}$ & $1.24\times10^{-4}$ & $3.57\times10^{-7}$ \\
64 & exact & 977 & $1.80\times10^{-8}$ & $2.12\times10^{-5}$ & $2.00\times10^{-8}$ & $2.95\times10^{-5}$ & $2.00\times10^{-8}$ \\
128 & exact & 4002 & $1.09\times10^{-9}$ & $5.18\times10^{-6}$ & $1.21\times10^{-9}$ & $7.18\times10^{-6}$ & $1.21\times10^{-9}$ \\
\midrule
16 & $10^{4}$ & 54 & $1.36\times10^{-5}$ & $3.87\times10^{-4}$ & $1.40\times10^{-5}$ & $5.42\times10^{-4}$ & $6.56\times10^{-5}$ \\
16 & $10^{5}$ & 54 & $8.43\times10^{-6}$ & $3.87\times10^{-4}$ & $8.95\times10^{-6}$ & $5.43\times10^{-4}$ & $1.38\times10^{-5}$ \\
16 & $10^{6}$ & 54 & $7.94\times10^{-6}$ & $3.87\times10^{-4}$ & $8.50\times10^{-6}$ & $5.43\times10^{-4}$ & $8.87\times10^{-6}$ \\
32 & $10^{4}$ & 233 & $2.31\times10^{-6}$ & $8.90\times10^{-5}$ & $2.32\times10^{-6}$ & $1.24\times10^{-4}$ & $7.35\times10^{-5}$ \\
32 & $10^{5}$ & 233 & $7.91\times10^{-7}$ & $8.93\times10^{-5}$ & $7.96\times10^{-7}$ & $1.24\times10^{-4}$ & $8.45\times10^{-6}$ \\
32 & $10^{6}$ & 233 & $3.95\times10^{-7}$ & $8.91\times10^{-5}$ & $4.23\times10^{-7}$ & $1.24\times10^{-4}$ & $1.25\times10^{-6}$ \\
\bottomrule
\end{tabular}
\captionbelow{Cubic B-spline dictionary with $N$ functions for the model of \cref{ex:theta_growth}, $\beta=1/2$, $T=1$, and $\dt=T/n$ with $\dt\|\cA_N\|\le1.5$: errors $\|\cP\,\Proj_Nf-\cK^Tf\|_{\cH}$ for the exponential, SDMD-1, SDMD-R, and SDMD-2 with the pointwise evaluation of $\cA^2\psi$ and with the energy form, for exact matrices and for empirical matrices from $m$ i.i.d.\ uniform samples (means over $50$ repetitions). Fitted slopes in $N$ over $N=32$--$128$ for exact matrices: $-4.10$, $-2.05$, $-4.10$, $-2.05$, $-4.10$.}
\label{tab:spline}
\end{table}

\FloatBarrier
\section{Conclusions}\label{sec:conclusion}
Among the second-order SDMD schemes, the extrapolated scheme SDMD-R should be used. It needs only the two matrices of SDMD-1 and a dictionary in $\dom(\cA)$, and its error against the Koopman semigroup is bounded by the Galerkin error of the dictionary plus a time-stepping term of order $\dt^2$ and a data term of order $m^{-1/2}$. The Taylor scheme SDMD-2 needs a dictionary in $\dom(\cA^2)$ and a third matrix, and its distance from the Galerkin semigroup, the common limit of the schemes as $\dt\to0$, is in general of first order. For reversible periodic diffusions, trigonometric dictionaries, and smooth observables the two schemes reach a given accuracy with the same cutoff, samples, and steps, so the cost of the matrices decides. In the experiments with trigonometric dictionaries their errors against the Koopman semigroup differ by at most ten percent, and for a spline dictionary with two derivatives only SDMD-R is second-order accurate with pointwise evaluations. The finite-sample bounds assume independent samples of the invariant law and a dictionary that is fixed or independent of the samples; dictionaries learned from the same samples and spectral conclusions for the generator are left for future work.

\appendix
\section{Two remarks on SDMD-1 and SDMD-2}\label{app:remarks}

\subsection{Relation to the SDMD paper}\label{app:relation}
The scheme of Xu et al.~\cite[Eq.~(10)]{SDMD2025} is the empirical step $\mathbf I+\dt\,\widehat{\mathbf G}_m^{-1}\widehat{\mathbf H}_m$ of \eqref{eq:empirical_matrices}. Since $\widehat{\mathbf G}_m^{-1}\widehat{\mathbf H}_m$ is the generator EDMD matrix of Klus et al.~\cite{KlusEtAl2020}, the eigenvectors of the two methods agree, and an eigenvalue $\lambda$ of the generator matrix corresponds to the eigenvalue $1+\dt\lambda$ of the step, so the conversion of Xu et al.~\cite[Sec.~III.B]{SDMD2025} returns the generator EDMD eigenvalues exactly. The analysis of \cref{sec:fixedN,sec:fourier} sharpens four points of Xu et al.~\cite{SDMD2025}. The large-sample limit of Eq.~(16) there is the projected Koopman step $\cK_N^{\dt}$, with matrix $\mathbf G^{-1}[\ip{\psi_i}{\cK^{\dt}\psi_j}]_{i,j=1}^N$. When $\cV_N\subset\dom(\cA^3)$, \eqref{eq:dynkin} gives $\cK_N^{\dt}=\cQ^{[1]}_{N,\dt}+\tfrac12\dt^2\cB_N+O(\dt^3)$ with $\cB_N$ from \eqref{eq:sdmd2}, so the zero-lag convergence of Theorem~IV.8 there holds with rate $O(\dt)$, and the time-truncation error does not vanish as $m\to\infty$ (related small-lag effects are noted in \cite[Sec.~4.3]{BramburgerFantuzzi2024} and \cite[Sec.~4.2]{Mauroy2021}). The concentration bound of Theorem~IV.3 there is replaced by the residual form of \cref{lem:residual,prop:perturbation}, which incorporates the regularization $\gamma$ used there and has a Bernstein-type rate; samples from a numerical invariant law are not covered (\cref{rem:data_variants}). Iterating the one-step map, as suggested in Remark~IV.16 there, converges at fixed $N$ to $e^{T\cA_N}$ as $\dt\to0$ and not to $\Proj_N\cK^T|_{\cV_N}$; for the projected Koopman step the leading correction is $\dt\,\cL_{N,T}$ (\cref{rem:reading}). The dictionary-limit argument of Theorem~IV.14 there needs a bound on $\|e^{t\cA_N}\|$ uniform in $N$, which \eqref{eq:AN_contraction} provides.

\subsection{An energy form for the third matrix}\label{app:energy}
If $\cA$ is self-adjoint, then $\mathbf J_{ij}=\ip{\cA\psi_i}{\cA\psi_j}$ can be estimated by $\widehat{\mathbf J}^{\mathrm{en}}_m\coloneqq m^{-1}(\cA\boldsymbol\Psi)^*(\cA\boldsymbol\Psi)\succeq0$ from the array $\cA\boldsymbol\Psi$ alone, in analogy with the reversible form of generator EDMD \cite{KlusEtAl2020}. This removes the evaluations of $\cA^2\psi$ from \cref{tab:cost} but not the $mN^2$ products of the third matrix. The form also defines $\cB_N$ for dictionaries in $\dom(\cA)$ only, through $\ip{u}{\cB_Nv}\coloneqq\ip{\cA u}{\cA v}$ for $u,v\in\cV_N$; with this definition $\cD_N=(\cE_N^\perp)^*\cE_N^\perp$ as in \cref{prop:residual_structure}, and \cref{lem:onestep} and the statements of \cref{thm:main} on SDMD-2 remain valid. For generators that are not self-adjoint the corresponding form involves the adjoint of $\cA$ in $L^2(\rho)$ and therefore the invariant density. Since $\widehat{\mathbf J}^{\mathrm{en}}_m$ is not of residual form, its error bound follows from the resolvent identity and carries the summand bound $F_1^2$ in place of $F_0F_s$ and an additional Gram term; in the Fourier setting of \cref{sec:fourier} both bounds scale like $K^{d+4}$ in dimension $d$. For the spline dictionary of \cref{subsec:num_spline}, on which the form is needed, the sampling error of $\widehat{\mathbf J}^{\mathrm{en}}_m$ dominates the data term of SDMD-2.

\begin{acknowledgements}
The author thanks Kathrin Hellmuth for helpful discussions during a visit to Caltech.
\end{acknowledgements}

\section*{Declarations}

\paragraph{Funding} This work was supported by JST CREST, grant numbers JPMJCR24Q1 and JPMJCR24Q6, including the AIP Challenge Program.

\paragraph{Competing interests} The author has no competing interests to declare that are relevant to the content of this article.

\paragraph{Data availability} The numerical experiments use only synthetic samples generated by the code listed under Code availability. No other data sets were generated or analyzed.

\paragraph{Code availability} The code that reproduces every table and figure is available at \url{https://github.com/TalkingDoll/sdmd_analysis}.

\paragraph{Use of generative AI} Generative AI tools were used to improve the writing. The author reviewed and verified all content and takes full responsibility for it.

\bibliographystyle{spmpsci}
\bibliography{reference}

@article{SDMD2025,
  author  = {Xu, Y. and Shao, K. and Ishikawa, I. and Hashimoto, Y. and Logothetis, N. and Shen, Z.},
  title   = {A data-driven framework for {Koopman} semigroup estimation in stochastic dynamical systems},
  journal = {Chaos: An Interdisciplinary Journal of Nonlinear Science},
  volume  = {35},
  number  = {10},
  pages   = {103123},
  year    = {2025},
  doi     = {10.1063/5.0283640}
}

@article{Koopman1931,
  author  = {Koopman, B. O.},
  title   = {Hamiltonian systems and transformation in {Hilbert} space},
  journal = {Proceedings of the National Academy of Sciences of the United States of America},
  volume  = {17},
  number  = {5},
  pages   = {315--318},
  year    = {1931},
  doi     = {10.1073/pnas.17.5.315}
}

@article{Mezic2005,
  author  = {Mezi\'{c}, I.},
  title   = {Spectral properties of dynamical systems, model reduction and decompositions},
  journal = {Nonlinear Dynamics},
  volume  = {41},
  number  = {1--3},
  pages   = {309--325},
  year    = {2005},
  doi     = {10.1007/s11071-005-2824-x}
}

@article{Schmid2010,
  author  = {Schmid, P. J.},
  title   = {Dynamic mode decomposition of numerical and experimental data},
  journal = {Journal of Fluid Mechanics},
  volume  = {656},
  pages   = {5--28},
  year    = {2010},
  doi     = {10.1017/S0022112010001217}
}

@article{Tu2014,
  author  = {Tu, J. H. and Rowley, C. W. and Luchtenburg, D. M. and Brunton, S. L. and Kutz, J. N.},
  title   = {On dynamic mode decomposition: Theory and applications},
  journal = {Journal of Computational Dynamics},
  volume  = {1},
  number  = {2},
  pages   = {391--421},
  year    = {2014},
  doi     = {10.3934/jcd.2014.1.391}
}

@book{Kutz2016,
  author    = {Kutz, J. N. and Brunton, S. L. and Brunton, B. W. and Proctor, J. L.},
  title     = {Dynamic Mode Decomposition: Data-Driven Modeling of Complex Systems},
  publisher = {SIAM},
  year      = {2016},
  doi       = {10.1137/1.9781611974508},
  isbn      = {978-1-61197-449-2}
}

@article{Williams2015,
  author  = {Williams, M. O. and Kevrekidis, I. G. and Rowley, C. W.},
  title   = {A data-driven approximation of the {Koopman} operator: Extending dynamic mode decomposition},
  journal = {Journal of Nonlinear Science},
  volume  = {25},
  number  = {6},
  pages   = {1307--1346},
  year    = {2015},
  doi     = {10.1007/s00332-015-9258-5}
}

@article{KordaMezic2018,
  author  = {Korda, M. and Mezi\'{c}, I.},
  title   = {On convergence of extended dynamic mode decomposition to the {Koopman} operator},
  journal = {Journal of Nonlinear Science},
  volume  = {28},
  number  = {2},
  pages   = {687--710},
  year    = {2018},
  doi     = {10.1007/s00332-017-9423-0}
}

@article{vcrnjaric2020koopman,
  author  = {{\v{C}}rnjari\'{c}-{\v{Z}}ic, N. and Ma\'{c}e\v{s}i\'{c}, S. and Mezi\'{c}, I.},
  title   = {{Koopman} operator spectrum for random dynamical systems},
  journal = {Journal of Nonlinear Science},
  volume  = {30},
  number  = {5},
  pages   = {2007--2056},
  year    = {2020},
  doi     = {10.1007/s00332-019-09582-z}
}

@article{wanner2022robust,
  author  = {Wanner, M. and Mezi\'{c}, I.},
  title   = {Robust approximation of the stochastic {Koopman} operator},
  journal = {SIAM Journal on Applied Dynamical Systems},
  volume  = {21},
  number  = {3},
  pages   = {1930--1951},
  year    = {2022},
  doi     = {10.1137/21M1414425}
}

@article{KlusEtAl2020,
  author  = {Klus, S. and N\"{u}ske, F. and Peitz, S. and Niemann, J.-H. and Clementi, C. and Sch\"{u}tte, C.},
  title   = {Data-driven approximation of the {Koopman} generator: Model reduction, system identification, and control},
  journal = {Physica D: Nonlinear Phenomena},
  volume  = {406},
  pages   = {132416},
  year    = {2020},
  doi     = {10.1016/j.physd.2020.132416}
}

@article{nuske2023finite,
  author  = {N\"{u}ske, F. and Peitz, S. and Philipp, F. and Schaller, M. and Worthmann, K.},
  title   = {Finite-data error bounds for {Koopman}-based prediction and control},
  journal = {Journal of Nonlinear Science},
  volume  = {33},
  number  = {1},
  pages   = {14},
  year    = {2023},
  doi     = {10.1007/s00332-022-09862-1}
}

@article{llamazares2024data,
  author  = {Llamazares-Elias, L. and Llamazares-Elias, S. and Latz, J. and Klus, S.},
  title   = {Data-driven approximation of {Koopman} operators and generators: Convergence rates and error bounds},
  journal = {Journal of Computational Dynamics},
  volume  = {15},
  pages   = {163--206},
  year    = {2026},
  doi     = {10.3934/jcd.2026012}
}

@article{philipp2026variance,
  author  = {Philipp, F. M. and Schaller, M. and Boshoff, S. and Peitz, S. and N\"{u}ske, F. and Worthmann, K.},
  title   = {Variance representations and convergence rates for data-driven approximations of {Koopman} operators},
  journal = {Physica D: Nonlinear Phenomena},
  volume  = {492},
  pages   = {135223},
  year    = {2026},
  doi     = {10.1016/j.physd.2026.135223}
}

@article{BramburgerFantuzzi2024,
  author  = {Bramburger, J. J. and Fantuzzi, G.},
  title   = {Auxiliary functions as {Koopman} observables: Data-driven analysis of dynamical systems via polynomial optimization},
  journal = {Journal of Nonlinear Science},
  volume  = {34},
  number  = {1},
  pages   = {8},
  year    = {2024},
  doi     = {10.1007/s00332-023-09990-2}
}

@article{Mauroy2021,
  author  = {Mauroy, A.},
  title   = {{Koopman} operator framework for spectral analysis and identification of infinite-dimensional systems},
  journal = {Mathematics},
  volume  = {9},
  number  = {19},
  pages   = {2495},
  year    = {2021},
  doi     = {10.3390/math9192495}
}

@book{KloedenPlaten1992,
  author    = {Kloeden, P. E. and Platen, E.},
  title     = {Numerical Solution of Stochastic Differential Equations},
  publisher = {Springer},
  series    = {Applications of Mathematics},
  volume    = {23},
  year      = {1992},
  doi       = {10.1007/978-3-662-12616-5},
  isbn      = {978-3-642-08107-1}
}

@article{Tropp2012,
  author  = {Tropp, J. A.},
  title   = {User-friendly tail bounds for sums of random matrices},
  journal = {Foundations of Computational Mathematics},
  volume  = {12},
  number  = {4},
  pages   = {389--434},
  year    = {2012},
  doi     = {10.1007/s10208-011-9099-z}
}

@article{talay1990invariant,
  author  = {Talay, D.},
  title   = {Second-order discretization schemes of stochastic differential systems for the computation of the invariant law},
  journal = {Stochastics and Stochastic Reports},
  volume  = {29},
  number  = {1},
  pages   = {13--36},
  year    = {1990},
  doi     = {10.1080/17442509008833606}
}

@article{mattingly2002ergodicity,
  author  = {Mattingly, J. C. and Stuart, A. M. and Higham, D. J.},
  title   = {Ergodicity for {SDEs} and approximations: Locally {Lipschitz} vector fields and degenerate noise},
  journal = {Stochastic Processes and their Applications},
  volume  = {101},
  number  = {2},
  pages   = {185--232},
  year    = {2002},
  doi     = {10.1016/S0304-4149(02)00150-3}
}

@article{paulin2015concentration,
  author  = {Paulin, D.},
  title   = {Concentration inequalities for {Markov} chains by {Marton} couplings and spectral methods},
  journal = {Electronic Journal of Probability},
  volume  = {20},
  number  = {79},
  pages   = {1--32},
  year    = {2015},
  doi     = {10.1214/EJP.v20-4039}
}

@book{Kato1995,
  author    = {Kato, T.},
  title     = {Perturbation Theory for Linear Operators},
  publisher = {Springer},
  series    = {Classics in Mathematics},
  edition   = {Second},
  year      = {1995},
  doi       = {10.1007/978-3-642-66282-9},
  isbn      = {978-3-540-58661-6}
}

@article{ito1998trotter,
  author  = {Ito, K. and Kappel, F.},
  title   = {The {Trotter--Kato} theorem and approximation of {PDEs}},
  journal = {Mathematics of Computation},
  volume  = {67},
  number  = {221},
  pages   = {21--44},
  year    = {1998},
  doi     = {10.1090/S0025-5718-98-00915-6}
}

@article{Higham2005,
  author  = {Higham, N. J.},
  title   = {The scaling and squaring method for the matrix exponential revisited},
  journal = {SIAM Journal on Matrix Analysis and Applications},
  volume  = {26},
  number  = {4},
  pages   = {1179--1193},
  year    = {2005},
  doi     = {10.1137/04061101X}
}

@article{Richardson1911,
  author  = {Richardson, L. F.},
  title   = {The approximate arithmetical solution by finite differences of physical problems involving differential equations, with an application to the stresses in a masonry dam},
  journal = {Philosophical Transactions of the Royal Society of London. Series A},
  volume  = {210},
  pages   = {307--357},
  year    = {1911},
  doi     = {10.1098/rsta.1911.0009}
}

@book{HairerNorsettWanner1993,
  author    = {Hairer, E. and N{\o}rsett, S. P. and Wanner, G.},
  title     = {Solving Ordinary Differential Equations {I}: Nonstiff Problems},
  publisher = {Springer},
  series    = {Springer Series in Computational Mathematics},
  volume    = {8},
  edition   = {Second},
  year      = {1993},
  doi       = {10.1007/978-3-540-78862-1},
  isbn      = {978-3-540-56670-0}
}

@misc{bartel2026extendeddynamicmodedecomposition,
  author        = {Bartel, F. and Ritter, S. and Schaller, M. and Worthmann, K.},
  title         = {Extended dynamic mode decomposition with {Fourier} dictionaries: {Error} bounds and fast implementation},
  year          = {2026},
  eprint        = {2609.20144},
  archiveprefix = {arXiv},
  primaryclass  = {math.DS},
  howpublished  = {{arXiv}:2609.20144},
  doi           = {10.48550/arXiv.2609.20144}
}

\end{document}